\documentclass[aos,MSNbibl,nameyear,dvips]{arximspdf}

\doi{10.1214/13-AOS1184} 
\volume{42}
\issue{1}
\pubyear{2014}
\firstpage{87}
\lastpage{114}

\makeatletter
\def\cal{\mathcal}
\newcommand{\eqref}[1]{(\ref{#1})}
\newcommand{\field}[1]{\mathbb{#1}}
\newcommand{\R}{\field{R}}
\newcommand{\p}{\field{P}}
\newcommand{\Z}{\field{Z}}
\newcommand{\E}{\field{E}}
\newcommand{\FF}{\mathcal{F}}
\newcommand{\PP}{\mathcal{P}}

\def\argmin{\mathop{\operatorname{argmin}}}
\newtheorem{theorem}{Theorem}
\newproclaim{example}{Example}
\newproclaim{remark}{Remark}

\newproclaim{definition}{Definition}
\newtheorem{corol}{Corollary}

\newtheorem{proposition}{Proposition}

\def\lf{\lfloor}
\def\rf{\rfloor}
\def\x{\mathbf{x}}
\def\ZZ{\mathbf{Z}}
\def\var{\operatorname{Var}}
\def\cov{\operatorname{Cov}}
\def\dmd{\diamondsuit}
\makeatother

\begin{document}
\begin{frontmatter}

\title{Inference of weighted $V$-statistics for nonstationary time series
and its applications}
\runtitle{Weighted $V$-statistics for nonstationary time series}

\begin{aug}
\author[A]{\fnms{Zhou} \snm{Zhou}\corref{}\thanksref{t2}\ead[label=e1]{zhou@utstat.toronto.edu}}
\thankstext{t2}{Supported by NSERC of Canada.}
\runauthor{Z. Zhou}
\affiliation{University of Toronto}
\address[A]{Department of Statistical Sciences\\
University of Toronto\\
100 St. George St.\\
Toronto Ontario\\
M5S3G3 Canada\\
\printead{e1}} 
\end{aug}

\received{\smonth{3} \syear{2013}}
\revised{\smonth{10} \syear{2013}}

%
\begin{abstract}
We investigate the behavior of Fourier transforms for a wide class of
nonstationary nonlinear processes. Asymptotic central and noncentral
limit theorems are established for a class of nondegenerate and
degenerate weighted $V$-statistics through the angle of Fourier analysis.
The established theory for $V$-statistics provides a unified treatment
for many important time and spectral domain problems in the analysis of
nonstationary time series, ranging from nonparametric estimation to
the inference of periodograms and spectral densities.
\end{abstract}

%

\end{frontmatter}

\section{Introduction}
\label{sec::intr}
Consider the following weighted $V$-statistics:
%
\begin{equation}
\label{eq:ustat} V_n=\sum_{k=1}^n
\sum_{j=1}^nW_n(t_k,t_j)H(X_k,X_j),
\end{equation}
where $\{X_k\}_{k=1}^n$ is a nonstationary time series, $t_k=k/n$, the
kernel $H(\cdot,\cdot)$ and the weights $W_n(\cdot,\cdot)$ are
symmetric and Borel measurable functions. Many important time and
spectral domain problems in the analysis of nonstationary time series
boil down to the investigation of weighted $V$-statistics in the form
of \eqref{eq:ustat}. For instance, in various situations one may be
interested in estimating parameter functions
\[
\theta(t)=\int_{\R}\int_{\R}H(u,v)
\,dF(t,u)\,dF(t,v)
\]
over time $t$, where $F(t,\cdot)$ is the marginal distribution of $\{
X_j\}$ at time $t$. In this case many nonparametric estimators of
$\theta(t)$ are asymptotically equivalent to the weighted $V$-statistics
\[
V_n=\sum_{k=1}^n\sum
_{j=1}^nK\bigl((t_k-t)/b_n,(t_j-t)/b_n
\bigr)H(X_k,X_j),
\]
where $K(\cdot,\cdot)$ is two-dimensional kernel function, and $b_n$ is
a bandwidth that restricts the estimation in a neighborhood of $t$.
Additionally, after the nonparametric fitting one may want to specify
or test whether the parameter function is of certain parametric forms.
In this case many ${\cal L}^2$ distance based test statistics are of
the form~\eqref{eq:ustat} with degenerate kernels; that is, kernels
$H(\cdot,\cdot)$ such that $\E H[X_k,x]= 0$ for every $k$ and $x$. See,
for instance, the ${\cal L}^2$ distance based quantile specification
test in \citet{Zho10}. Furthermore, in spectral analysis of the
nonstationary time series $\{X_k\}$, both the periodogram
\[
I_n(\lambda)=\frac{1}{2\pi n}\Biggl|\sum_{j=1}^nX_j
\exp(ij\lambda)\Biggr|^2, \qquad 0\le\lambda\le\pi,
\]
and the classic smoothed periodogram estimate of the spectral density
\[
\tilde{f}_n(\lambda)=\int_{\R}\frac{1}{m}K
\biggl(\frac{u}{m} \biggr)I_n(\lambda+2\pi u/n) \,du
\]
are in the form of \eqref{eq:ustat} with kernel $H(x,y)=xy$. Here
$i=\sqrt{-1}$ stands for the imaginary unit, $K(\cdot)$ is a kernel
function and $m=m_n$ is a window size satisfying $m\rightarrow\infty$
with $m/n\rightarrow0$. Note that in every example above, $W_n(\cdot
,\cdot)$ involves a tuning parameter (either bandwidth or window size)
which varies with the sample size $n$. Hence it is important to write
the weights as a function of $n$ in \eqref{eq:ustat}.

The purpose of the paper is to establish an asymptotic theory for
\eqref
{eq:ustat} through the angle of Fourier analysis. To illustrate the
main idea, suppose that $H^*(x,y)=H(x,y)/(L(x)L(y))$ is absolutely
integrable on $\R^2$ for some function $L(\cdot)$. Then under mild
conditions $H^*(\cdot,\cdot)$ admits the Fourier representation
$H^*(x,y)=\int_{\R^2}g(u,v)e^{i(xu+yv)} \,du \,dv$. Consequently $V_n$ can
be written as
%
\begin{equation}
\label{eq:fourier_repre} V_n=\int_{\R^2}g(x,y)\sum
_{k,j=1}^nW_n(t_k,t_j)
\beta_k(x)\beta _j(y) \,dx \,dy,
\end{equation}
where $\beta_k(x)=L(X_k)\exp(ixX_k)$. Note that in \eqref
{eq:fourier_repre} the complex structure of $V_n$ is reduced to a
process of quadratic forms $\{\sum_{k,j}W_n(t_k,t_j)\beta_k(x)\beta
_j(y)\}_{x,y\in\R}$ through the aid of Fourier transformation. The
multiplicative structure of the quadratic forms makes in-depth
asymptotic investigations possible for a wide class of nonstationary
time series. Furthermore, the continuous structure of $\beta_k(x)$ in
$x$ makes a stochastic equi-continuity and continuous mapping argument
possible which is shown to be powerful compared to the discrete
spectral decomposition methods used in the literature. With the aid of
the above structural simplifications, in this paper we are able to
establish a uniform approximation scheme of $\{\sum_{k,j}W_n(t_k,t_j)\beta_k(x)\beta_j(y)\}_{x,y\in\R}$ by a process of
Gaussian quadratic forms. As a consequence we establish a unified
asymptotic theory for a class of nondegenerate and degenerate weighted
$V$-statistics with reflexible weight functions. Both central and
noncentral limit theorems are developed for a wide class of
nonstationary time series with both smoothly and abruptly changing
data generating mechanisms over time. The established theory can be
applied to many problems in the study of nonstationary time series,
including topics such as nonparametric estimation and specification,
periodogram and spectral density inferences discussed above.

In \eqref{eq:ustat}, if the summation is taken over indices $1\le
k\neq
j\le n$, then the statistics are commonly called weighted
$U$-statistics in the literature. For many problems that are of
practical importance, statistical theory for the $U$- and $V$-statistics
can be established with essentially the same techniques. Since the
seminal papers of \citet{vMi47} and \citet{Hoe48}, the analysis
of $V$- and $U$-statistics has attracted much attention in the
statistics and probability literature. It seems that most efforts have
been put in un-weighted $U$- or $V$-statistics with $W_n(\cdot, \cdot
)\equiv1$ for stationary data. See, for instance, \citet{Yos76},
\citet{DehTaq89}, Huskova and Janssen (\citeyear{HusJan93}), \citet{DehWen10}, \citet{Leu12} and
Beutner and Z\"{a}hle (\citeyear{BeuZah12,BeuZah13}) for
various approaches for nondegenerate and degenerate un-weighted $U$-
and $V$-statistics. We also refer to the monographs of \citet{Den85},
\citet{Lee90} and \citet{Deh06} for more references. As we observe from
the examples in the beginning of this \hyperref[sec::intr]{Introduction}, it is important to
consider weighted $V$-statistics with sample size dependent weights in
the study of nonstationary time series.

There are a small number of papers discussing weighted $V$- or
$U$-statistics in the literature. See, for instance, \citet{deJ87},
\citet{ONeRed93}, \citet{Maj94}
and \citet{RifUtz00} for weighted $U$- and $V$-statistics of
independent data and \citet{HsiWu04} for weighted nondegenerate
$U$-statistics of stationary time series. For most of the above
discussions, the weights are not allowed to be sample size dependent.
Exceptions include \citet{deJ87} who discovered a very deep result
that $\theta_{n,1}\rightarrow0$ implies asymptotic normality of a very
wide class of weighted degenerate $V$-statistics for independent data,
where $\theta_{n,1}$ is the eigenvalue of the matrix $\{
W_{n}(t_j,t_k)/\sqrt{\sum_{u,v=1}^nW^2_{n}(t_u,t_v)}\}_{j,k=1}^n$ with
the maximum absolute value. There the proof heavily depended on the
martingale structure of degenerate $V$-statistics of independent data and
is hard to generalize to the time series setting. In this paper, from a
Fourier analysis angle, we generalize the result of \citet{deJ87} and
show that, for many temporally dependent processes, $\theta
_{n,1}\rightarrow0$ implies asymptotic normality of $V_n$ for a class
of degenerate kernels and weight functions.

Quadratic forms are special cases of \eqref{eq:ustat} with $H(x,y)=xy$.
There are many papers in the literature devoted to the analysis of such
statistics. See, for instance, \citet{deWVen73}, \citet{FoxTaq87}, G\"{o}tze and Tikhomirov (\citeyear{GotTik99}), \citet{GaoAnh00} and
Bhansali, Giraitis and
Kokoszka (\citeyear{BhaGirKok07}), among others. It seems that most of the results
are on independent or stationary data. Exceptions include \citet{LeeSub} who recently studied asymptotic normality of a class of
quadratic forms with banded weight matrix for $\alpha$-mixing
nonstationary time series. For independent data, G\"{o}tze and
Tikhomirov (\citeyear{GotTik99}), among others, established deep theoretical results
indicating that distributions of generic quadratic forms can be
approximated by those of corresponding Gaussian quadratic forms. In
this paper, we generalize this type of result and show that the laws of
a wide class of quadratic forms for nonstationary time series can be
well approximated by the distributions of corresponding quadratic forms
of independent Gaussian random variables. Consequently, central and
noncentral limit theorems are established for the latter class of
quadratic forms for nonstationary time series.

The rest of the paper is organized as follows. In Section~\ref{sec:pre}
we shall introduce the class of absolutely convergent Fourier
transformations and the piece-wise locally stationary time series
models used in this paper. Sections~\ref{sec:nondegen} and \ref
{sec:degen} establish the asymptotic theory for nondegenerate and
degenerate $V$-statistics, respectively. Theory for quadratic forms will
be covered in Section~\ref{sec:degen} as a special case of degenerate
$V$-statistics. In Section~\ref{sec:app}, we shall apply our theory to
the problems of nonparametric estimation as well as spectral analysis
of nonstationary processes. Several examples will be discussed in
detail. Finally, the theoretical results are proved in Section~\ref{sec:proof}.



\section{Preliminaries}\label{sec:pre}
We first introduce some notation. For a vector $\mathbf{
v} = (v_1, v_2,  \ldots, v_p) \in\R^p$, let $|\mathbf{v}| =
(\sum_{i=1}^p v_i^2)^{1/2}$. Let $i=\sqrt{-1}$ be the imaginary unit.
For a complex number $z=x+yi \in\field{C}$, write $|z|=\sqrt {x^2+y^2}$. For $q>0$, denote by $L^q(\R^p)$ the collection of
functions $f \dvtx \R^p\rightarrow\field{C}$ such that $\int_{\R
^p}|f(\mathbf{
x})|^q \,d\mathbf{x}<\infty$. For a function $f\in L^1(\R^p)$, denote by
$\hat{f}$ its Fourier transform, that is, $\hat{f}(\mathbf{v})=\int_{\R
^p}f(\x)e^{-i\langle\x,\mathbf{v}\rangle} \,d\x$. For a Borel set $A$
in $\R
^p$, denote by $\mathfrak{B}(A)$ the collection of all Borel sets in
$A$. For a random vector
$\mathbf{V}$, write $\mathbf{V} \in{\cal L}^q$ ($q > 0$) if $\|\mathbf{
V}\|_q := [\E(|\mathbf{V}|^q) ]^{1/q} < \infty$ and $\|\mathbf{V}\| =
\|\mathbf{V}\|_2$. Denote by $\Rightarrow$ the weak convergence. The
symbol $C$ denotes a generic finite constant which may vary from place
to place.

\subsection{Absolutely convergent Fourier transforms}

Following the classic notation, a function $f \dvtx \R^2\rightarrow\R$ is
said to belong to $W_0(\R^2)$ if there exists a function~$g$, $\R
^2\rightarrow\field{C}$, such that $g\in L^1(\R^2)$ and
%
\begin{equation}
\label{eq:w0} f(x,y)=\int_{\R^2} g(t,s)e^{itx+isy} \,dt \,ds.
\end{equation}
The class $W_0(\R^2)$ is called the Wiener ring or Wiener algebra,
naming after Norbert Wiener for his fundamental contributions in the
study of absolutely convergent Fourier integrals. Due to its importance
in various problems, the possibility to represent a function as an
absolutely convergent
Fourier integral has been intensively studied in mathematics, and
various sufficient conditions are continuously being discovered in
recent years. For a relatively comprehensive survey, we refer to
Liflyand, Samko and
  Trigub (\citeyear{LifSamTri12}).

For the study of $V$-statistics of dependent data, we need to further
define the subclass $W_0^{\delta}(\R^2)$ of $W_0(\R^2)$ as follows. We
call $f\in W_0^{\delta}(\R^2)$ for some $\delta\ge0$ if $f\in
W_0(\R
^2)$ with representation \eqref{eq:w0} and
%
\begin{equation}
\label{eq:w0delta} \int_{\R^2}\bigl|(t,s)\bigr|^{\delta}\bigl|g(t,s)\bigr| \,dt\,ds<
\infty,
\end{equation}
where $0^0:=1$. Note that $W_0^0(\R^2)=W_0(\R^2)$.

In order for $f\in W_0(\R^2)$, it is necessary that $f$ is uniformly
continuous on $\R^2$ and vanishes at $\infty$. However, the above
conditions are generally not sufficient. Generally speaking, the
Fourier transform of a smoother function tends to have a lighter tail
at infinity. The latter fact implies that stronger smoothness
conditions are needed to insure $f\in W_0(\R^2)$.
%
\begin{proposition}\label{prop:ac}
Suppose that a symmetric function $f\in L^1(\R^2)$. \textup{(i)} [Liflyand, Samko and
  Trigub (\citeyear{LifSamTri12}), Theorem 10.11]. Assume that $f$ is uniformly continuous on
$\R^2$ and vanishes at $\infty$. Let $f$ and its partial derivative
$(\partial/\partial x)f$ be locally absolutely continuous on $(\R
\setminus\{0\})^2$ in each variable. Further let each partial
derivative $(\partial/\partial x)f$ and $[\partial^2/(\partial
x\,\partial y)]f$ exist and belong to $L^p(\R^2)$ for some $p\in(1,
\infty
)$. Then $f\in W_0(\R^2)$. \textup{(ii)} If $f$ satisfies
%
\begin{eqnarray}
\label{eq:1}&& \biggl|\frac{\partial^2 }{\partial x^2}f(x+\epsilon,y)-\frac
{\partial^2
}{\partial x^2}f(x,y)\biggr |+ \biggl|
\frac{\partial^2 }{\partial
x^2}f(x,y+\epsilon)-\frac{\partial^2 }{\partial x^2}f(x,y) \biggr|
\nonumber\\
&&\quad{}+ \biggl|
\frac{\partial^2 }{\partial x\, \partial y}f(x+\epsilon ,y)-\frac
{\partial^2 }{\partial x\, \partial y}f(x,y) \biggr|\\
&&\qquad\le C|\epsilon
|^{\gamma} k(x,y)\nonumber
\end{eqnarray}
for sufficiently small $\epsilon$, where $\gamma>0$ and $k\in L^1(\R
^2)$, then $f\in W_0^{\delta}(\R^2)$ for any $\delta\in[0,\gamma)$.
\end{proposition}

Proposition \ref{prop:ac} gives some easily checkable sufficient
conditions for $f\in W_0(\R^2)$ and $W_0^{\delta}(\R^2)$ based on the
smoothness of its partial derivatives. In the literature, numerous
other sufficient conditions based on various notions of smoothness or
variation are available; see, for instance, the review of Liflyand, Samko and
  Trigub (\citeyear{LifSamTri12}). We point out here that the conditions in Proposition \ref
{prop:ac} are not minimal sufficient. For instance, the function
$f=\exp
(-|\x|)\in W_0^{\delta}(\R^2)$ for any $\delta\in[0,1/2)$. But the
latter function is not differentiable on $\R^2$. Thanks to the fast
computation of Fourier transforms, in practice when facing specific
choice of $f$, condition \eqref{eq:w0delta} with $g=\hat{f}$ can also
be checked via numerical computation.

\subsection{Nonstationary time series models}
For an observed process\vspace*{1pt} $\{X_j\}_{j=1}^n$, consider the class of
nonstationary time series models of the form [Zhou (\citeyear{Zho13}), Definition 1]
%
\begin{equation}
\label{eq:std_cons} X_k=\sum_{j=0}^rG_j(t_k,
\FF_k)I_{(b_j,b_{j+1}]}(t_k),\qquad  k=1,2,\ldots,n,
\end{equation}
where $b_1<b_2<\cdots< b_r$ are $r$ unknown (but nonrandom) break
points with $b_0=0$, and $b_{r+1}=1$, $I_{\cdot}(\cdot)$ is the
indicator function, $G_j\dvtx (b_j,b_{j+1}]\times\R^\field{N}\rightarrow
\R
$ are $(\mathfrak{B}((b_j, b_{j+1}])\times\mathfrak{B}(\R)^\field{N}$,
$\mathfrak{B}(\R))$-measurable functions, $j=0,\ldots,r$, $\FF
_k=(\ldots
, \varepsilon_{k-1},\varepsilon_k)$ and $\varepsilon_k$'s are i.i.d.
random variables. Recall that $t_k=k/n$. We observe from \eqref
{eq:std_cons} that the data generating mechanism could break at the
points $b_k$, $k=1,\ldots,r$ and hence lead to structural breaks of $\{
X_j\}$ at the latter points. On the other hand, note that if
$G_k(t,\cdot)$ are smooth functions of $t$ for each $k$, then the data
generating mechanisms change smoothly between adjacent break points. As
a consequence $\{X_j\}$ is approximately stationary in any small
temporal interval between adjacent break points. By the above
discussion, we shall call the class of time series in the form of
\eqref
{eq:std_cons} piece-wise locally stationary (PLS) processes.

Time series in the form of \eqref{eq:std_cons} constitute a relatively
large class of nonstationary time series models which allow the data
generating mechanism to change flexibly over time. In particular, when
the number of the break points $r=0$, then \eqref{eq:std_cons} reduces
to the locally stationary time series models in \citet{ZhoWu09}. On
the other hand, if for each $k=0,1,\ldots,r$, $G_k(t,\FF_j)$ does not
depend on $t$, then \eqref{eq:std_cons} is a piece-wise stationary time
series which is studied, for instance, in Davis, Lee and
  Rodriguez-Yam (\citeyear{DavLeeRod06}), among
others. Process \eqref{eq:std_cons} can be viewed as a time-varying
nonlinear system with $\varepsilon_k$'s being the inputs and $X_k$'s
being the outputs. The functions $G_k(t,\ldots)$ can be viewed as
time-varying filters of the system. From this point of view, we adapt
the following dependence measures of $\{X_k\}$ in Zhou (\citeyear{Zho13}):
%
\begin{equation}
\label{eq:pdm} \delta(j, p) = \sup_{t}\max
_k \bigl\|G_k({t,\cal F}_k) -
G_k(t,{\cal F}_{k, k-j})\bigr \|_p,
\end{equation}
where ${\cal F}_{j, k}$ is a coupled version of ${\cal F}_j$
with $\varepsilon_k'$ in ${\cal F}_j$ replaced by an i.i.d. copy
$\varepsilon_k'$; that is,
\[
{\cal F}_{j, k} = \bigl(\ldots, \varepsilon_{k-1},
\varepsilon_k', \varepsilon_{k+1}, \ldots,
\varepsilon_{j-1}, \varepsilon_j\bigr),
\]
and $\{\varepsilon'_j\}$ is an i.i.d. copy of $\{\varepsilon_j\}$.

We observe from \eqref{eq:pdm} that $\delta(j,p)$ measures the impact
of the system's inputs $j$-steps before on the current output of the
system. When $\delta_{j,p}$ decays fast to zero as~$j$ tends to
infinity, we have short memory of the series as the system tends to
fast ``forget about'' past inputs. We refer to Section~2 on page 728 of
Zhou (\citeyear{Zho13}) and Section~4 on pages 2706--2708 of \citet{ZhoWu09}
for more examples of linear and nonlinear nonstationary time series of
the form \eqref{eq:std_cons} and detailed calculations of the
dependence measures \eqref{eq:pdm}.

\section{Nondegenerate $V$-statistics}\label{sec:nondegen}
Suppose that, for some function $L(\cdot)$,
$H^*(s,t):=H(s,t)/[L(s)L(t)]\in W_0(\R^2)$ and $\max_j\|L(X_j)\|
<\infty
$. Then $V_n$ defined in \eqref{eq:ustat} admits the representation
\eqref{eq:fourier_repre} with $\int_{\R^2}|g(x,y)| \,dx \,dy<\infty$.
Now define
\[
H_j(x)=\E\bigl[H(x,X_j)\bigr]=\E\bigl[H(X_j,x)
\bigr]=\int_{\R^2}g(t,s)L(x)\exp(itx)\E \bigl[\beta
_j(s)\bigr] \,dt\,ds
\]
and $\gamma_j(x)=\beta_j(x)-\E[\beta_j(x)]$. Then elementary
calculations using the Hoeffding's decomposition show that $V_n$ can be
decomposed as
%
\begin{equation}
\label{eq:decomp} V_n-\E V_n=2N_n+D_n-
\E[D_n],
\end{equation}
where
\begin{eqnarray*}
N_n&=&\sum_{k=1}^n\sum
_{j=1}^nW_n(t_k,t_j)
\bigl\{H_j(X_k)-\E\bigl[H_j(X_k)
\bigr]\bigr\} ,
\\
D_n&=&\int_{\R^2}g(x,y)\sum
_{k,j}W_n(t_k,t_j)
\gamma_k(x)\gamma _j(y) \,dx \,dy.
\end{eqnarray*}
Here $N_n$ and $D_n$ are the nondegenerate and degenerate part of
$V_n$, respectively.

To investigate the limiting behavior of $N_n$, we need the following conditions:
\begin{longlist}[(A1)]
\item[(A1)] For some $\eta\in(0,1]$, $H^*(t,s):=H(t,s)/[L(t)L(s)]\in
W_0^{\eta}(\R^2)$ for some function $L(\cdot)$.
\item[(A2)] The function $L$ is differentiable with derivative $L'$.
$\max_{a\le t\le b}|L'(t)|\le C(|L'(a)|+|L'(b)|+1)$ for all $a$, $b$,
where $C$ is a finite constant independent of $a$ and $b$. $\max_j\|
L(X_j)\|_{4+2\epsilon}<\infty$ and $\max_j\|L'(X_j)\|_{4+2\epsilon
}<\infty$ for some $\epsilon>0$.
%
\item[(A3)] Define $W_{j,\cdot}=\sum_{r=1}^n|W_n(t_j,t_r)|$,
$j=1,2,\ldots,n$. Let $W_{j,\cdot}=0$ for $j>n$. Let $W^{(n)}=\sum_{j=1}^nW^2_{j,\cdot}$. For sequences $l_n$, $m_n$ and $s_n=l_n+m_n$, define
\begin{eqnarray}
A_j=\sum_{k=1}^{l_n}W^2_{s_n(j-1)+k,\cdot}\quad
\mbox{and}\quad a_j=\sum_{k=l_n+1}^{s_n}W^2_{s_n(j-1)+k,\cdot},\nonumber\\
\eqntext{j=1,2,\ldots,\lceil n/s_n\rceil.}
\end{eqnarray}
Assume that there exist sequences
$m_n/\log n\rightarrow\infty$ with $l_n/n\rightarrow0$ and
$m_n/l_n\rightarrow0$, such that
\[
\sum_{j}a_j/W^{(n)}
\rightarrow0 \quad\mbox{and}\quad \max_jA_j/W^{(n)}
\rightarrow0.
\]
\item[(A4)] $\phi_n:=\var(N_n)/W^{(n)}\ge c$ for some $c>0$ and
sufficiently large $n$.

\item[(A5)] The dependence measures $\delta_X(k,4+2\epsilon)=O(\rho^k)$
for some $\rho\in[0,1)$ and $\epsilon>0$.
\end{longlist}

A few comments on the regularity conditions are in order. The role of
$L(t)$ in (A1) is to lighten the tail of $H(s,t)$ and hence make it
absolutely integrable on $\R^2$. Some typical choices of $L(t)$ are
$(1+t^2)^p$ for kernels $H(\cdot,\cdot)$ with algebraically increasing
tails and $\exp[(1+t^2)^p]$ for kernels $H(\cdot,\cdot)$ with
exponentially increasing tails. Since the function $L(t)$ controls the
tail behavior of the kernel $H(\cdot,\cdot)$, condition (A1)
essentially poses some smoothness requirement on $H(\cdot,\cdot)$. In
practice, Proposition \ref{prop:ac} or direct numerical computations
can be used to check $H^*(t,s)\in W_0^{\delta}(\R^2)$.

Condition (A2) posts some moment restrictions on $L(X_j)$ and
$L'(X_j)$. The requirement $\max_{a\le t\le b}|L'(t)|\le
C(|L'(a)|+|L'(b)|+1)$ is mild and in particular it is always satisfied
when $|L'(\cdot)|$ is piece-wise monotone or piece-wise convex with
finite many pieces. Consequently, the latter inequality holds for
functions $L(t)=(1+t^2)^p$ or $\exp[(1+t^2)^p]$ listed above. Condition
(A3) requires that $\{1,2,\ldots, n\}$ can be divided into big blocks
and small blocks such that the sum of squares of $W_{j,\cdot}$ in the
small blocks is negligible and $W_{j,\cdot}$ in the big blocks satisfy
a Lindeberg type condition. Condition (A3) can be checked easily in
practice. Examples \ref{ex:1}--\ref{ex:3} below verify (A3) for some
frequently used weight functions. By the proof of Theorem \ref
{thm:clt_nondegen} in Section~\ref{sec:proof}, $\var(N_n)=O(W^{(n)})$.
Condition (A4) avoids degenerate kernels with $N_n\equiv0$ and other
degenerate cases with non-Gaussian limits. For instance, if
$H(x,y)=x+y$, $W_n(x,y)\equiv1/n$ and $X_j=Y_{j}-Y_{j-1}$, where $Y_j$
is a weakly stationary time series, then we have that $N_n=Y_n-Y_0$
fails to be asymptotically normal. Note in this case $W^{(n)}=n$ and
$\var(N_n)/W^{(n)}\rightarrow0$. Finally (A5) requires that the
dependence measures of $\{X_j\}$ decay exponentially fast to zero. The
theoretical results of the paper can be derived with $\delta_X(k,p)$
decaying polynomially fast to zero at the expense of much more
complicated conditions and proofs. For presentational simplicity and
clarity we assume exponentially decaying dependence measures throughout
the paper.

\begin{theorem}\label{thm:clt_nondegen}
Under conditions \textup{(A1)--(A5)}, we have
\[
N_n/\sqrt{\var(N_n)}\Rightarrow N(0,1).
\]
\end{theorem}
Theorem \ref{thm:clt_nondegen} establishes a central limit result for
the nondegenerate part of~$V_n$. Due to nonstationarity, $H_j(x)=\E
[H(x,X_j)]$ is $j$-dependent. Consequently asymptotic investigations
for $N_n$ are much more difficult than the stationary case. Thanks to
the structural simplification by the Fourier transformation, we are
able to control the dependence structure of $N_n$ and establish the
above central limit results. By Theorem \ref{thm:clt_nondegen}, if the
degenerate part $D_n$ is asymptotically negligible compared to $N_n$,
then CLT for $V_n$ can be derived. To this end, the following Theorem~\ref{thm:order_degen} is crucial.

\begin{longlist}[(A6)]
\item[(A6)] Assume that there exist functions $f_n$ and a $p\in
(0,\infty
)$, such that
\[
\bigl|W_n(t_l,t_m)-W_n(t_k,t_j)\bigr|
\le f_n(t_k,t_j)\bigl|(l,m)-(k,j)\bigr|^p
\]
for all integers $k,j,l,m\in\{1,2,\ldots,n\}$.
\end{longlist}

Condition (A6) requires that, for each fixed $n$, the weight function
$W_n(x,y)$ grows algebraically fast at any $(x,y)$. (A6) is mild and is
satisfied by most frequently used weight functions. In particular, if
$W_n(\cdot, \cdot)$ is Lipschitz continuous in the sense that
$|W_n(t,s)-W_n(t_1,s)|\le q_n|t_1-t|$ for every $t,t_1$ and $s$, then
$f_n(\cdot,\cdot)$ can be chosen as $Cq_n/n$.

\begin{theorem}\label{thm:order_degen}
Assume that conditions \textup{(A2)} and \textup{(A5)} hold with $4+2\epsilon$
therein
replaced by $8$. Further assume \textup{(A1)} and \textup{(A6)}. Then
\[
\|D_n-\E D_n\|^2=O(W_{(n)}+
\Delta_n) \qquad\mbox{where } W_{(n)}=\sum
_{k=1}^n\sum_{j=1}^nW^2_n(t_k,t_j)
\]
and $\Delta_n=\sum_{k=1}^n\sum_{j=1}^n|W_n(t_k,t_j)|f_n(t_k,t_j)$.
\end{theorem}
Theorem \ref{thm:order_degen} investigates the order of $D_n$. Observe
that if $X_j$'s are independent with $H(x,y)=xy$, then simple
calculations yield $\|D_n-\E D_n\|^2=O(W_{(n)})$. For many important
weight functions (see, e.g., Examples \ref{ex:1}--\ref{ex:3}
below) $\Delta_n=O(W_{(n)})$ and the order established in Theorem \ref
{thm:order_degen} is sharp. For un-weighted $U$-statistics of stationary
mixing data, deep theoretical results on the order of the degenerate
part of $U$-statistics were obtained in \citet{Yos76} and \citet{DehWen10}, among others. Our Theorem \ref{thm:order_degen} extends
the latter results to weighted $V$-statistics of nonstationary processes
from a Fourier analysis angle.

\begin{corol}\label{coro:clt_nondegen}
Assume that conditions \textup{(A2)} and \textup{(A5)} hold with $4+2\epsilon$ therein
replaced by $8$. Further assume that \textup{(A1), (A3), (A4)} and \textup{(A6)} hold and
$W^{(n)}/[W_{(n)}+\Delta_n]\rightarrow
\infty$. Then we have $(V_n-\E V_n)/\sqrt{\var(V_n)}\Rightarrow N(0,1)$.
\end{corol}
Corollary \ref{coro:clt_nondegen} is an immediate consequence of
Theorems \ref{thm:clt_nondegen} and \ref{thm:order_degen}. Corollary
\ref{coro:clt_nondegen} establishes a CLT for $V_n$. It is clear from
the definitions of $W^{(n)}$ and $W_{(n)}$ that $W^{(n)}\ge W_{(n)}$.
Examples \ref{ex:1}--\ref{ex:3} below verify
$W^{(n)}/[W_{(n)}+\Delta
_n]\rightarrow\infty$ and condition~(A3) for some frequently used
weight functions.

\begin{example}\label{ex:1}
Consider the case where $W_n(t_j,t_k)=f(t_j,t_k)/n$ for some symmetric
function $f$ on $[0,1]\times[0,1]$ such that $|f(x_1,y)-f(x_2,y)|\le
C|x_1-x_2|$ for all $x_1,x_2$ and $y$ in $[0,1]$. Note that the classic
un-weighted $V$-statistics are contained in this case with $f(\cdot
,\cdot
)\equiv1$. Elementary calculations yield $W_{j,\cdot}=f(t_j,\cdot
)+O(1/n)$, where $f(t,\cdot)=\int_{0}^1|f(t,x)| \,dx$, $W^{(n)}= n\int_{0}^1(\int_{0}^1|f(x,\break y)| \,dy)^2 \,dx+O(1)$, $W_{(n)}=O(1)$ and $\Delta
_n=O(1/n)$. Hence $W^{(n)}/[W_{(n)}+\Delta_n]\rightarrow
\infty$ provided that $f$ is not a constant zero function.
Additionally, it is elementary to check that~(A3) is satisfied for
every sequence $m_n/\log n\rightarrow\infty$, $l_n/n\rightarrow0$ and
$m_n/l_n\rightarrow0$ provided that $f$ is not always $0$.
\end{example}

\begin{example}\label{ex:2}
In this example we investigate weight functions in the form
$W_n(t_j,t_k)=g((t_j-t)/b_n,(t_k-t)/b_n)/(nb_n)$, where $t\in[0,1]$,
$g(\cdot,\cdot)$ is a continuously differentiable function and
$b_n\rightarrow0$ with $nb_n\rightarrow\infty$. This type of weights
may appear in nonparametric estimation of nonstationary time series.
Assume that $g(\cdot,\cdot)$ is absolutely integrable and its first-order partial derivatives are bounded on $\R^2$. Then elementary
calculations show that $W^{(n)}=(nb_n)\int_{\R}(\int_{\R}|g(x,y)|
\,dy)^2 \,dx+O(1)$, $W_{(n)}=O(1)$ and $\Delta_n=O(1/(nb_n))$ for any
$t\in
(0,1)$. Therefore $W^{(n)}/[W_{(n)}+\Delta_n]\rightarrow
\infty$ provided that $g$ is not always zero and (A3) is satisfied for
every sequence $m_n/\log n\rightarrow\infty$, $l_n/(nb_n)\rightarrow
0$ and $m_n/l_n\rightarrow0$ provided that $g$ is not always $0$.
Similar results hold for $t=0$ or $1$.
\end{example}

\begin{example}\label{ex:3}
Consider the class of weights $W_n(t_j,t_k)=\sqrt {m_n}h(|t_j-t_k|m_n)/ n$, where $h(\cdot)$ is a continuously
differentiable and nonconstant function on $[0,\infty)$ and
$m_n\rightarrow\infty$ with $m_n/n\rightarrow0$. This type of weight
functions may appear, for instance, in nonparametric specification
tests and spectral analysis of $\{X_j\}$. Further assume that $h$ is
absolutely integrable on $[0,\infty)$ with bounded derivatives. After
some simple algebra, we have $W^{(n)}=\frac{n}{m_n}[4(\int_{0}^\infty
|h(x)| \,dx)^2+o(1)]$, $W_{(n)}=O(1)$ and $\Delta_n=O(m_n/n)$. Therefore
$W^{(n)}/[W_{(n)}+\Delta_n]\rightarrow
\infty$ and (A3) is satisfied for every sequence $m_n/\log
n\rightarrow
\infty$, $l_n/n\rightarrow0$ and $m_n/l_n\rightarrow0$.
\end{example}

\section{Degenerate $V$-statistics}\label{sec:degen}
Without loss of generality, throughout Section~\ref{sec:degen} we
assume $c\le\sum_{k=1}^n\sum_{j=1}^nW^2_n(t_k,t_j)\le C$ for some
constants $0<c\le C<\infty$. In this section we shall investigate the
class of degenerate $V$-statistics for which $\E[H(x,X_j)]=0$ for every
$j$ and $x$. Then it is clear that $N_n=0$ in \eqref{eq:decomp} and
$V_n=D_n$. Before we state the theoretical results, we need to post the
following regularity conditions:

\begin{longlist}[(A7)]
\item[(A7)] Define
\[
{\cal V}_n=\sum_{j=1}^{n-1}
\Biggl[\sum_{k=1}^n\bigl(W_n(t_j,t_k)-W_n(t_{j+1},t_{k})
\bigr)^2 \Biggr]^{1/2}+ \Biggl[\sum
_{k=1}^nW^2_n(1,t_k)
\Biggr]^{1/2}.
\]
Assume that $n^{1/4}\log^2 n{\cal V}_n=o(1)$. Further assume that for
some $\delta>0$,
%
\begin{equation}
\label{eq:A7} \sum_{j=1}^n \biggl[\sum
_{|k-j|\le\log^{1+\delta}
n}\bigl(W_n(t_j,t_j)-W_n(t_{j},t_{k})
\bigr)^2 \biggr]^{1/2}=o(1).
\end{equation}

\item[(A8)] If $|m-l|=O(\log n)$, then $\sum_{k=m}^l[\sum_{|j-k|\le
\log
^{1+\delta}n}W^2_n(t_k,t_j)]^{1/2}=o(1)$ for some $\delta>0$.

\item[(A9)] For $s=(s_1,s_2,\ldots,s_m)^{\top}\in\R^m$ and $t\in
[0,1]$, define
\begin{eqnarray*}
&&\beta_{k,j}(t,s)=L\bigl(G_k(t,\FF_j)\bigr)
\bigl(\cos\bigl(s_1G_k(t,\FF_j)\bigr),\sin
\bigl(s_1G_k(t,\FF_j)\bigr),\ldots,\\
&&\hspace*{134pt}\cos
\bigl(s_mG_k(t,\FF_j)\bigr),\sin
\bigl(s_mG_k(t,\FF_j)\bigr)
\bigr)^{\top}.
\end{eqnarray*}
Let $\beta^*_{k,j}(t,s)=(L(G_k(t,\FF_j)),\beta^\top
_{k,j}(t,s))^\top$.
For $k=0,1,\ldots,r$, define the long-run covariances
$\Sigma_k(t,s)=\sum_{j=-\infty}^\infty\cov[\beta_{k,0}(t,s),\beta
_{k,j}(t,s)]$ and $\Sigma^*_k(t,s)=\break \sum_{j=-\infty}^\infty\cov
[\beta
^*_{k,0}(t,s),\beta^*_{k,j}(t,s)]$. Assume that for all $m\in\field
{N}$ and all $s\in\R^m$ with $s_1<s_2<\cdots<s_m$ and $s_j\neq0$,
$\Sigma_k(t,s)$ is positive definite for $t\in[b_k,b_{k+1}]$,
$k=0,1,\ldots,r$ if $L(x)\equiv C$. Replace $\Sigma_k(t,s)$ by
$\Sigma
^*_k(t,s)$ in the above assumption for all other functions $L(x)$.

\item[(A10)] $\|G_k(t,\FF_0)-G_k(s,\FF_0)\|_4\le C|t-s|$ for $t,s\in
[b_{k},b_{k+1}]$ and $k=0,1,\ldots,r$.
\end{longlist}

Condition (A7) posts some restrictions on the smoothness of the weight
function $W_n(\cdot,\cdot)$. Condition (A8) is a mild technical
condition. In particular, elementary calculations show that (A7) and
(A8) are satisfied by weight functions in Example~\ref{ex:1}. We have
(A7) and (A8) hold if $b_n\gg n^{-1/2}\log^4 n$ under the additional
assumption that $\int_{\R}(\int_{\R}(\frac{\partial
g(x,y)}{\partial
y})^2 \,dy)^{1/2} \,dx<\infty$ in Example \ref{ex:2}. For weights
considered in Example \ref{ex:3}, we have (A7) and (A8) hold when
$m_n\ll n^{1/4}/\log^2 n$ under the extra assumption that $\int_{0}^\infty(h'(x))^2 \,dx<\infty$.

By definition, $\Sigma_k(t,s)$ and $\Sigma^*_k(t,s)$ are the spectral
density matrices of the time series $\{\beta_{k,j}(t,s)\}
_{j=0}^{\infty
}$ and $\{\beta^*_{k,j}(t,s)\}_{j=0}^{\infty}$ at frequency 0,
respectively. Hence it is clear that the latter spectral density
matrices are positive semi-definite. Condition~(A9) is mild, and it
requires that the latter spectral density matrices are nonsingular.
Finally, condition (A10) means that the data generating mechanism
$G_k(t,\cdot)$ changes smoothly between adjacent break points. The
following theorem investigates the asymptotic behavior of degenerate
weighted $V$-statistics:

\begin{theorem}\label{thm:degen}
Write $L^*(x)=xL(x)$. Let condition \textup{(A2$'$)} be condition \textup{(A2)} when $L(x)$
therein is replaced by $L^*(x)$. Assume that conditions \textup{(A2)}, \textup{(A2$'$)} and
\textup{(A5)} hold with $4+2\epsilon$ therein replaced by $8+4\epsilon$ for some
$\epsilon>0$. Assume \textup{(A1)} and \textup{(A7)--(A10)}. Then there exist constants
$\alpha_{n,1},\alpha_{n,2},\ldots $ with $\sum_{k=1}^\infty\alpha
^2_{n,k}=O(1)$ and i.i.d. standard normal random variables
$Z_1,Z_2,\ldots,$ such that for any bounded and continuous function
$h(\cdot)$
%
\begin{equation}
\label{eq:degen} \Biggl|\E\bigl[h(D_n-\E D_n)\bigr]-\E\Biggl[h
\Biggl(\sum_{j=1}^\infty\alpha _{n,j}
\bigl(Z_j^2-1\bigr)\Biggr)\Biggr] \Biggr|\rightarrow0.
\end{equation}
\end{theorem}
Let $\Gamma(t,s)=(\Gamma_{1}(t,s),\Gamma_{2}(t,s))^\top$ be a centered
two-dimensional Gaussian process defined on $[0,1]\times\R$ with the
covariance function
\[
\cov\bigl[\Gamma(t_1,s_1),\Gamma(t_2,s_2)
\bigr]=\int_0^{\min(t_1,t_2)}\Xi _{\zeta
(t)}
\bigl(t,(s_1,s_2)^{\top}\bigr) \,dt,
\]
where $\zeta(t)=k$ if $b_k<t\le b_{k+1}$, $k=0,1,\ldots,r$, $\zeta
(0)=0$ and
\[
\Xi_k\bigl(t,(s_1,s_2)^{\top}\bigr)=
\sum_{j=-\infty}^\infty\cov\bigl[\beta
_{k,0}(t,s_1),\beta_{k,j}(t,s_2)
\bigr].
\]
Define the complex-valued Gaussian process $\Gamma^*(t,s)=\Gamma
_{1}(t,s)+i\Gamma_{2}(t,s)$ and let
%
\begin{equation}
\label{eq:40} \Gamma_n^*(t,s)=\sqrt{n}\biggl[\Gamma^*(t,s)-\Gamma^*
\biggl(t-\frac
{1}{n},s\biggr)\biggr], \qquad t\ge1/n.
\end{equation}
Let $\varpi_n(x,y)=\sum_{k,j=1}^nW_n(t_k,t_j)\Gamma_n^*(t_k,x)\Gamma
_n^*(t_j,y)$. Then by the classic Gaussian process theory [see, for
instance, \citet{Kuo75}, Chapter~1.2], the real part of $\int_{\R
^2}g(x,y)\varpi_n(x,y) \,dx\,dy$ is a quadratic form of i.i.d. Gaussian
random variables $Z_1,Z_2,\ldots.$ The coefficients $\alpha_{n,j}$,
$j=1,2,\ldots,$ in Theorem \ref{thm:degen} correspond to the eigenvalues
of the latter Gaussian quadratic form. Theorem \ref{thm:degen}
establishes a general asymptotic result for degenerate $V$-statistics
with smooth weight functions. Define the L\'{e}vy--Prokhorov metric
\[
\pi(\mu,\nu)=\inf\bigl\{\epsilon>0| \mu(A)\le\nu\bigl(A^{\epsilon
}\bigr)+
\epsilon, \nu (A)\le\mu\bigl(A^{\epsilon}\bigr)+\epsilon\mbox{ for every Borel
set }A\bigr\},
\]
where $A^{\epsilon}$ is the $\epsilon$-neighborhood of $A$. Then
\eqref
{eq:degen} is equivalent to
\[
\pi\Biggl(\operatorname{law}(D_n-\E D_n), \operatorname{law}\Biggl(\sum
_{j=1}^\infty\alpha _{n,j}
\bigl(Z_j^2-1\bigr)\Biggr)\Biggr)\rightarrow0.
\]
In other words, the distribution of $D_n-\E D_n$ can be well
approximated asymptotically by that of a weighted sum of i.i.d.
centered $\chi^2(1)$ random variables. An important observation from
\eqref{eq:degen} is that if $\max_{1\le j<\infty}|\alpha
_{n,j}|\rightarrow0$, then $\pi(\operatorname{law}(\sum_{j=0}^\infty\alpha
_{n,j}(Z_j^2-1)), \operatorname{law}(N(0,2\sum_{j=0}^\infty\alpha
_{n,j}^2)))\rightarrow0$. Consequently $D_n-\E D_n$ is asymptotically
normal. On the other hand, if $\alpha_{n,j}\rightarrow\alpha_j$
uniformly in $j$ as $n\rightarrow\infty$, then $\sum_{j=1}^\infty
\alpha
_{n,j}(Z_j^2-1))\rightarrow\sum_{j=1}^\infty\alpha_{j}(Z_j^2-1))$ and
therefore $D_n-\E D_n$ converges to a mixture of i.i.d. centered $\chi
^2(1)$ random variables. In the following, Corollaries \ref
{coro:clt_degen} to \ref{coro:nonclt-degen2} explore the above
discussions in detail.

\begin{corol}\label{coro:clt_degen}
Let $\theta_{n,1},\theta_{n,2},\ldots,\theta_{n,n}$ be the eigenvalues
of the matrix $\{W_n(t_j,t_k)\}_{j,k=1,\ldots,n}$ with $|\theta
_{n,1}|\ge|\theta_{n,2}|\ge|\theta_{n,n}|$. Assume that $\theta
_{n,1}\rightarrow0$.
Then under the conditions of Theorem \ref{thm:degen}, we have for any
bounded and continuous function~$h(\cdot)$
%
\begin{equation}
\label{eq:CLT_degen} \bigl|\E h[D_n-\E D_n] -\E h\bigl\{ N\bigl(0,
\var[D_n]\bigr)\bigr\}\bigr|\rightarrow0.
\end{equation}
\end{corol}
Note that by Lemma 4 in Zhou (\citeyear{Zho13}), $\var[D_n]=O(1)$. 
Corollary \ref{coro:clt_degen} asserts that $|\theta
_{n,1}|\rightarrow
0$ implies $\max_{1\le j<\infty}|\alpha_{n,j}|\rightarrow0$ and hence
the asymptotical normality of~$D_n$. In the literature, \,\citet{deJ87}
derived that $\theta_{n,1}\rightarrow0$ implies asymptotic normality
of a very wide class of weighted degenerate $V$-statistics of independent
data based on very deep martingale techniques. The martingale arguments
depended heavily on the independence assumption and are hard to
generalize to dependent data. From a Fourier analysis point of view,
Corollary \ref{coro:clt_degen} generalizes the latter result to a class
of weighted degenerate $V$-statistics of nonstationary time series with
smooth weights. A particular case of this type is Example~\ref{ex:3}
where the weight matrix is Toeplitz. From standard Toeplitz matrix
theory we have $|\theta_{n,1}|\le\sum_{j=1}^n\sqrt {m_n}|h(t_jm_n)|/n=O(1/\sqrt{m_n})=o(1)$, and hence \eqref
{eq:CLT_degen} holds for this type of weight matrices.

\begin{corol}\label{coro:nonclt-degen1}
Suppose that \textup{(a)}: $W_n(t,s)=Q_1(t,s)/n$ and $Q_1(t,s)$ satisfies \eqref
{eq:1} with $\R^2$ therein replaced by $[0,1]\times[0,1]$; or \textup{(b)}:
$W_n(t,s)=\sum_{j=1}^\infty a_{j}\times f_{1,j}(t)f_{2,j}(s)/n$, where $\sum_{j=1}^\infty|a_j|<\infty$ and $f_{1,j}(\cdot)$ and $f_{2,j}(\cdot)$
are continuous functions defined on $[0,1]$. Then under the conditions
of Theorem \ref{thm:degen}, there exist constants $\alpha_1,\alpha
_2,\ldots$ with $\sum_{j=1}^\infty\alpha_j^2<\infty$ and i.i.d.
standard normal random variables $Z_1,Z_2,\ldots,$ such that
\[
D_n-\E D_n\Rightarrow\sum_{j=1}^\infty
\alpha_j\bigl(Z_j^2-1\bigr).
\]
\end{corol}

\begin{corol}\label{coro:nonclt-degen2}
Suppose \textup{(a)}: $W_n(t,s)=Q_2((t-a)/b_n,(s-b)/b_n)/(nb_n)$ for some
$a,b\in[0,1]$, where $Q_2(t,s)$ has support $[-1,1]\times[-1,1]$ and
satisfies \eqref{eq:1} with $\R^2$ therein replaced by $[-1,1]\times
[-1,1]$ and $b_n\rightarrow0$; or \textup{(b)}:~$W_n(t,s)=\sum_{j=1}^\infty
a_jg_{1,j}((t-a)/b_n)g_{2,j}((s-b)/b_n)/(nb_n)$ for some $a,b\in
[0,1]$, where $\sum_{j=1}^\infty|a_j|<\infty$, $b_n\rightarrow0$ and
$g_{1,j}(\cdot)$ and $g_{2,j}(\cdot)$ are continuous functions on $\R$
with support $[-1,1]$. Then under the conditions of Theorem \ref
{thm:degen}, there exist constants $\alpha_1,\alpha_2,\ldots$ with
$\sum_{j=1}^\infty\alpha_j^2<\infty$ and\vadjust{\goodbreak} i.i.d. standard normal random
variables $Z_1,Z_2,\ldots,$ such that
\[
D_n-\E D_n\Rightarrow\sum_{j=1}^\infty
\alpha_j\bigl(Z_j^2-1\bigr).
\]
\end{corol}
Corollaries \ref{coro:nonclt-degen1} and \ref{coro:nonclt-degen2} are
proved in the online supplement of the paper, Zhou (\citeyear{Zho}). Corollaries
\ref{coro:nonclt-degen1} and \ref{coro:nonclt-degen2} establish that
$D_n$ converges to a mixture of i.i.d. centered $\chi^2(1)$ random
variables for four classes of smooth weight functions which are
absolutely integrable on $\R^2$. Note that the classic un-weighted
$V$-statistics belong to cases (a) and (b) in Corollary \ref
{coro:nonclt-degen1}. In the literature, \citet{Leu12}, among others,
derived asymptotic distributions of un-weighted $U$-statistics for
stationary time series. Corollary \ref{coro:nonclt-degen1} generalizes
the latter results to a class of weighted $V$-statistics of
nonstationary data. Weight functions in Corollary \ref
{coro:nonclt-degen2} may appear, for instance, in nonparametric
estimation of nonstationary time series.

Quadratic forms of a centered nonstationary process $\{X_j\}$ are of
the form
%
\begin{equation}
\label{eq:quad} Q_n=\sum_{j=1}^n
\sum_{k=1}^nW_n(t_j,t_k)X_jX_k.
\end{equation}
Clearly $Q_n$ is a special case of the degenerate $V$-statistics with
$H(x,y)=xy$. Hence the theory established above applies to this class
of statistics. However, due to the special multiplicative structure,
the asymptotic theory for $Q_n$ can be established with weaker
conditions. The following proposition follows from the corresponding
proofs of Theorem \ref{thm:degen} and Corollaries \ref{coro:clt_degen}
to \ref{coro:nonclt-degen2}.
%
\begin{proposition}\label{prop:quad}
Assume that conditions \textup{(A7), (A8)} and \textup{(A10)} hold and \textup{(A5)}
holds with
$4+2\epsilon$ therein replaced by $4$. Further assume that $\tilde
{\sigma}^2(k,t):=\sum_{j=-\infty}^\infty\cov[G_k(t,\FF
_0),G_k(t,\FF
_j)]>0$ for $k=0,1,\ldots, r$ and $t\in[b_k,b_{k+1}]$. Write $\tilde
{\sigma}^2(t)=\tilde{\sigma}^2(\zeta(t),t)$. Then we have that
conclusions of Theorem \ref{thm:degen} and Corollaries \ref
{coro:clt_degen} to \ref{coro:nonclt-degen2} hold with $D_n$ therein
replaced by $Q_n$.
\end{proposition}
By the proof of Theorem \ref{thm:degen}, on a possibly richer
probability space, there exist i.i.d. standard normal random variables
$Z_1,Z_2,\ldots, Z_n$, such that
%
\begin{equation}
\label{eq:quad_approx} Q_n-Q^o_n=o_\p(1),\qquad
\mbox{where }Q_n^o=\sum_{j,k=1}^nW_n(t_j,t_k)
\tilde {\sigma}(t_j)Z_j\tilde{\sigma}(t_k)Z_k.
\end{equation}
The above equation asserts that $Q_n$ can be well approximated a
quadratic form of independent Gaussian random variables. In the
literature, G\"{o}tze and Tikhomirov (\citeyear{GotTik99}), among others, established
deep theoretical results showing that distributions of quadratic forms
of independent data can be approximated by those of corresponding
Gaussian quadratic forms. In \eqref{eq:quad_approx}, we generalize
these results to a class of quadratic forms of nonstationary time
series with smooth weights.


\section{Applications}\label{sec:app}

\subsection{Nonparametric estimation of nonstationary time series}
Let $F(t,\cdot)$ be the marginal distribution function of $\{X_j\}$ at
time $t$; namely $F(t,\cdot)$ is the distribution of $G_{\zeta
(t)}(t,\FF
_0)$. Under various situations one is interested in estimating the quantity
%
\begin{equation}
\label{eq:estimation} \theta(t)=\int_{-\infty}^\infty\int
_{-\infty}^\infty H(x,y) \,d F(t, x)\,d F(t, y)
\end{equation}
for all $t\in[0,1]$. Here $H(\cdot,\cdot)$ is assumed to be a symmetric
function. For instance, if $H(x,y)=(x-y)^2/2$, then $\theta(t)$ is the
time-varying variance function of the process $\{X_j\}$. In the
statistics literature, enormous efforts have been put on nonparametric
estimation of parameter functions $\tau(t)$ in the form $\tau(t)=\int_{-\infty}^\infty M(x)\,d F(t, x)=E[M[G_{\zeta(t)}(t,\FF_0)]]$; see, for
instance, the monographs of \citet{FanGij96} and \citet{FanYao03} and the citations therein. Note that $\tau(t)$ is a special case
of \eqref{eq:estimation} with $H(x,y)=[M(x)+M(y)]/2$. On the other
hand, however, it seems that there are few results on nonparametric
inference of general parameter functions in the form of $\theta(t)$ in
\eqref{eq:estimation}. One of the major difficulties, especially in the
case of time series applications, lies in the lack of corresponding
theoretical results on weighted $V$-statistics for dependent data. Define
%
\begin{equation}
\label{eq:estimation_2d} \theta(t,s)=\int_{-\infty}^\infty\int
_{-\infty}^\infty H(x,y) \,d F(t, x)\,d F(s, y).
\end{equation}
Assume $\theta(t,s)$ is smooth at $(t^*,t^*)$ for some $t^*\in(0,1)$.
By the first-order local Taylor expansion of $\theta(t,s)$, $\theta
(t^*)$ can be estimated by $\hat{\theta}_{b_n}(t^*)$, where
%
\begin{eqnarray}
\label{eq:local_linear}&& \bigl(\hat{\theta}_{b_n}\bigl(t^*\bigr), \hat{
\eta}_1,\hat{\eta}_2\bigr)
\nonumber
\\
&&\qquad= \argmin_{(\eta_0, \eta_1 \eta_2)\in\R^3}
\sum_{j,k=1}^n \bigl(H(X_j,X_k)-
\eta_0
-\eta_1\bigl(t_j-t^*\bigr)\\
&&\hspace*{183pt}{}-
\eta_2\bigl(t_k-t^*\bigr) \bigr)^2
W_n(t_j,t_k).\nonumber
\end{eqnarray}
Here for presentational simplicity we assume that
\[
W_n(t_j,t_k)=K\bigl(\bigl(t_j-t^*
\bigr)/b_n\bigr)K\bigl(\bigl(t_k-t^*\bigr)/b_n
\bigr)/(nb_n),
\]
where $K\in{\cal K}$ and ${\cal K}$ is the collection of continuously
differentiable and symmetric density functions with support $[-1,1]$.
Furthermore, the bandwidth $b_n$ satisfies $b_n\rightarrow0$ with
$nb_n\rightarrow\infty$. Estimator \eqref{eq:local_linear} is an
extension of the classic local linear kernel methods [\citet{FanGij96}] to second-order parameter functions of the form \eqref
{eq:estimation}. Meanwhile, if higher-order Taylor expansions of
$\theta
(t,s)$ are used in~\eqref{eq:local_linear}, then one obtains local
polynomial estimations of $\theta(t^*)$.

It is easy to see that the asymptotic behavior of $\hat{\theta
}_{b_n}(t^*)$ is decided by that of the $V$-statistics $V_n=\sum_{j,k=1}^nH(X_j,X_k)W_n(t_j,t_k)$. The following proposition, which is
proved in Zhou (\citeyear{Zho13}), investigates the limiting distribution of $\hat
{\theta}_{b_n}(t^*)$:
%
\begin{proposition}\label{prop:local-linear}
Assume that conditions \textup{(A2)} and \textup{(A5)} hold with $4+2\epsilon$ therein
replaced by $8$. Further assume that \textup{(A1)} and \textup{(A4)} hold and $\theta
(t,s)$ is ${\mathcal C}^2$ in a neighborhood of $(t^*,t^*)$. Then under
the above assumptions of $K(\cdot)$ and $b_n$, we have
%
\begin{equation}
\label{eq:clt_local_linear} \frac{\sqrt{nb_n}}{\sqrt{4\phi_n\int_{-1}^1K^2(x) \,dx}}\bigl[\hat {\theta }_{b_n}\bigl(t^*
\bigr)-\theta\bigl(t^*\bigr)-B_n\bigl(t^*\bigr)\bigr]\Rightarrow
N(0,1),
\end{equation}
where $B_n(t^*)=b_n^2\frac{\partial^2 \theta(t^*,t^*)}{\partial
t^2}\int_{-1}^1x^2K(x) \,dx$.
\end{proposition}


\begin{example}[(Estimating the time-varying variance function)]
Consider the kernel $H(x,y)=(x-y)^2/2$. Then $\theta(t)$ in \eqref
{eq:estimation} equals $\var[G_{\zeta(t)}(t,\FF_0)]=\E[G_{\zeta
(t)}(t,\FF_0)-\E G_{\zeta(t)}(t,\FF_0)]^2$. In particular, $\theta
(t_i)=\var[X_i]$.

For this variance kernel $H$, we can choose $L(x)=(1+x^2)^2$ and assume
that $\E[X_i^{32}]<\infty$, $i=1,2,\ldots, n$. Then $L(x)$ satisfies
condition (A2) with $4+2\epsilon$ therein replaced by 8. By Proposition
\ref{prop:ac}, condition (A1) is satisfied with $\eta=1$. Furthermore,
condition (A3) is satisfied by Example \ref{ex:2} and the assumption that $K\in
{\cal K}$. Note that $H_j(x)=\E[H(x,X_j)]=(x^2-2x\E[X_j]+\E[X^2_j])/2$
does not always equal $0$. Hence the kernel is nondegenerate. Meanwhile,
\[
2\theta(t,s)=\E\bigl[G^2_{\zeta(t)}(t,\FF_0)\bigr]+
\E\bigl[G^2_{\zeta(s)}(s,\FF _0)\bigr]-2\E
\bigl[G_{\zeta(t)}(t,\FF_0)\bigr]\E\bigl[G_{\zeta(s)}(s,
\FF_0)\bigr].
\]
Assuming that $\mu(t):=\E[G_{\zeta(t)}(t,\FF_0)]$ and $v(t):=\E
[G^2_{\zeta(t)}(t,\FF_0)]$ are ${\cal C}^2$ in a neighborhood of $t^*$,
then $\theta(t,s)$ is ${\cal C}^2$ in a neighborhood of $(t^*,t^*)$.
Further assume that $t^*$ is not a break point of the time series and
condition (A10). By the local stationarity of $\{X_j\}$ in the
neighborhood of $t^*$, we have that $4\phi_n\rightarrow\sigma^2(t^*)$,
where
\begin{eqnarray*}
\sigma^2(t)&=&\sum_{j=-\infty}^\infty
\cov\bigl\{G^2_{\zeta(t)}(t,\FF _0)-2\mu
(t)G_{\zeta(t)}(t,\FF_0),G^2_{\zeta(t)}(t,
\FF_j)\\
&&\hspace*{160pt}{}-2\mu(t)G_{\zeta
(t)}(t,\FF_j)\bigr\}.
\end{eqnarray*}
Note that $\sigma^2(t)$ is the spectral density of the stationary
sequence $\{G^2_{\zeta(t)}(t,\break \FF_j)-2\mu(t)G_{\zeta(t)}(t,\FF_j)\}
_{j=-\infty}^\infty$ at frequency $0$. Hence condition (A4) holds
provided $\sigma^2(t^*)>0$. Finally, we have under the other regularity
assumptions of Proposition~\ref{prop:local-linear} that
\[
\frac{\sqrt{nb_n}}{\sqrt{\sigma^2(t^*)\int_{-1}^1K^2(x) \,dx}}\bigl[\hat {\theta}_{b_n}\bigl(t^*\bigr)-\theta
\bigl(t^*\bigr)-\tilde{B}_n\bigl(t^*\bigr)\bigr]\Rightarrow N(0,1),
\]
where $\tilde{B}_n(t^*)=\frac{b_n^2}{2}[v''(t^*)-2\mu(t^*)\mu
''(t^*)]\int_{-1}^1x^2K(x) \,dx$.
\end{example}

\subsection{Spectral analysis}
Consider a PLS time series $\{X_j\}$ defined in \eqref{eq:std_cons}.
Assume further that $\E[X_j]=0$ and $\|X_j\|<\infty$, $j=1,2,\ldots,n$.
Then we can define its spectral density at time $t$ as
\[
f(t,\lambda)=\frac{1}{2\pi}\sum_{k=-\infty}^\infty
\gamma (t,k)\cos (k\lambda), \qquad \lambda\in[0,2\pi],
\]
where $\gamma(t,k)=\cov[G_{\zeta(t)}(t,\FF_0),G_{\zeta(t)}(t,\FF_k)]$
is the $k$th-order auto covariance of $\{X_j\}$ at time $t$. Write the
classic periodogram of the series $\{X_j\}$
\[
I_n(\lambda)=\frac{1}{2\pi n}\bigl|S_n(\lambda)\bigr|^2,\qquad
\mbox{where }S_n(\lambda )=\sum_{j=1}^nX_j
\exp(ij\lambda),  0\le\lambda\le\pi.
\]
Consider also the classic smoothed periodogram estimate of the spectral density
\[
\tilde{f}_n(\lambda)=\int_{-m}^m
\frac{1}{m}K \biggl(\frac{u}{m} \biggr)I_n(\lambda+2\pi
u/n) \,du,
\]
where $K(\cdot)\in{\cal K}$ is an even function, and $m$ is a block
size satisfying $m\rightarrow\infty$ with $n/m\rightarrow\infty$. The
analysis of $\tilde{f}_n(\lambda)$ depends heavily on the theory of
quadratic forms for nonstationary processes. For strictly stationary
time series, the asymptotic behaviors of the periodogram and spectral
density estimates have been intensively studied in the literature. See,
for instance, \citet{Bri69}, \citet{Pri81}, \citet{Ros84},
\citet{Wal00} and \citet{ShaWu07} among others. On the other hand,
however, there are few corresponding results for nonstationary time
series. Exceptions include, among others, \citet{DwiSub11}
who studied the asymptotic behavior of the periodogram for short memory
locally stationary linear processes and Dette, Preuss and Vetter (\citeyear{DetPreVet11}) who studied
the behavior of the averaged spectral density estimates for locally
stationary Gaussian linear processes. In this section we shall
investigate the behaviors of $I_n(\lambda)$ and $\tilde{f}_n(\lambda)$
for linear and nonlinear PLS time series. The following is a key
theorem which establishes a Gaussian approximation result for Fourier
transforms of nonstationary time series. Theorem \ref
{thm:ftransform_normal_appr} could be of separate interest in spectral
analysis of nonstationary processes.

Let $\lambda=\lambda_n$ be a sequence of frequencies of interest. For
$1\le a <b\le n$, define $S^*_{a,b,\lambda}=\sum_{j=a}^bX_j\cos
(j\lambda
)$ and $S^o_{a,b,\lambda}=\sum_{j=a}^bX_j\sin(j\lambda)$. Write
$S^*_{a,\lambda}:= S^*_{1,a,\lambda}$ and $S^o_{a,\lambda}:=
S^o_{1,a,\lambda}$.
%
\begin{theorem}\label{thm:ftransform_normal_appr}
Assume that $\gamma(t,h)$ is ${\cal C}^p$ in $t$ on $(b_k,b_{k+1}]$ for
any $h$, $k=0,1,\ldots,r$, $p\ge1$. Further assume that conditions
\textup{(A5)} and \textup{(A10)} hold with $4+2\epsilon$ therein replaced
by 4 and $\inf_{t\in[0,1]}f(t,\lambda)>0$. \textup{(i)}: If $0<\lambda_*\le\lambda\le
\lambda
^*<\pi$ for some constants $\lambda_*$ and $\lambda^*$, then on a
possibly richer probability space, there exist i.i.d. two-dimensional
standard normal random vectors $(G_{1,1},G_{1,2})^\top, \ldots,
(G_{n,1},G_{n,2})^\top$, such that
\begin{eqnarray*}
&&\max_{1\le j\le n}\Biggl|\bigl(S^*_{j,\lambda},S^o_{j,\lambda}
\bigr)^\top-\sum_{k=1}^j\sqrt{
\pi}\bigl(f^{1/2}(t_k,\lambda)G_{k,1},f^{1/2}(t_k,
\lambda )G_{k,2}\bigr)^\top\Biggr|\\
&&\qquad=o_\p
\bigl(n^{p^*}\log^2 n\bigr),
\end{eqnarray*}
where $p^*=(p+4)/[4(p+2)]$. \textup{(ii):} If $\lambda=0$ or $\pi$, then on a
possibly richer probability space, there exist i.i.d. standard normal
random variables $G^*_1,\ldots,\break  G^*_n$, such that
\[
\max_{1\le j\le n}\Biggl|S^*_{j,\lambda}-\sqrt{2\pi}\sum
_{k=1}^jf^{1/2}(t_k,
\lambda)G^*_{k}\Biggr|=o_\p\bigl(n^{1/4}
\log^2 n\bigr).
\]
\end{theorem}
Based on Theorem \ref{thm:ftransform_normal_appr}, we have the
following corollary on the behavior of the periodogram for
nonstationary time series.
%
\begin{corol}\label{coro:periodogram}
Under the conditions of Theorem \ref{thm:ftransform_normal_appr}, we
have \textup{(i)}: if the frequency $\lambda$ satisfies $0<\lambda_*\le
\lambda
\le\lambda^*<\pi$, then
\[
I_n(\lambda)/\int_{0}^1f(t,
\lambda) \,dt\Rightarrow\operatorname{Exp}(1),
\]
where $\operatorname{Exp}(1)$ stands for the exponential distribution with mean 1.
\textup{(ii)}
If $\lambda=0$ or~$\pi$, then $I_n(\lambda)/\int_{0}^1f(t,\lambda
) \,dt$
converges in distribution to a $\chi^2(1)$ random variable.
\end{corol}

\begin{remark}
The condition $0<\lambda_*\le\lambda\le\lambda^*<\pi$ is important
for the validity of Theorem \ref{thm:ftransform_normal_appr} and
Corollary \ref{coro:periodogram}. By (ii) of Lemma 5 in Zhou (\citeyear{Zho13}),
we have if $1\le k\le C$ for some finite constant $C$, then
%
\begin{equation}
\label{eq:im_corr} \Biggl|\cov\bigl(S^*_{n,\lambda_k},S^o_{n,\lambda_k}
\bigr)-\pi\sum_{j=1}^n f(t_j,
\lambda _k)\sin(4\pi kt_j)\Biggr|=O\bigl(\log^2 n
\bigr),
\end{equation}
where $\lambda_k=2\pi k/n$. Therefore the real and imaginary parts of
$S_n(\lambda_k)$ are no longer uncorrelated, and the periodogram does
not converge to an $\operatorname{Exp}(1)$ distribution. Similar results hold for
frequencies near $\pi$. This is drastically different from the
stationary case where it is well known that the real and imaginary
parts of $S_n(\lambda_k)$ are asymptotically independent and
$I_n(\lambda_k)/f(\lambda_k)\Rightarrow\operatorname{Exp}(1)$. Indeed, note
that if $f(t,\lambda_k)$ does not change with $t$, then we have
$\sum_{j=1}^n f(t_j,\lambda_k)\sin(4\pi kt_j)=0$ in~\eqref{eq:im_corr}. Due\vadjust{\goodbreak}
to the time-varying nature of $f(t,\lambda_k)$, the behavior of Fourier
transforms near frequency 0 or $\pi$ is complicated for nonstationary
time series.
\end{remark}
The following proposition investigates the asymptotic behavior of
$\tilde{f}_n(\lambda)$ for PLS time series.
%
\begin{proposition}\label{prop:spectral-density}
Assume that $K\in{\cal K}$ is even. Then under the conditions of
Theorem 4 and the assumption that $m\rightarrow\infty$ with
$m/(n^{p/[4(p+2)]}\log^2 n)\rightarrow0$, we have \textup{(i)}: if $0<\lambda
_*\le\lambda\le\lambda^*<\pi$, then
\[
\sqrt{m}\bigl(\tilde{f}_n(\lambda)-\E\tilde{f}_n(
\lambda)\bigr)\Rightarrow N\biggl(0, \int_{-1}^1
\bigl[K(t)\bigr]^2 \,dt\int_{0}^1f^2(t,
\lambda) \,dt\biggr);
\]
and \textup{(ii)}: if $\lambda=0$ or $\pi$, then
\[
\sqrt{m}\bigl(\tilde{f}_n(\lambda)-\E\tilde{f}_n(
\lambda)\bigr)\Rightarrow N\biggl(0, 2\int_{-1}^1
\bigl[K(t)\bigr]^2 \,dt\int_{0}^1f^2(t,
\lambda) \,dt\biggr).
\]
\end{proposition}
Proposition \ref{prop:spectral-density}, which is proved in Zhou
(\citeyear{Zho13}), establishes the asymptotic normality of $\tilde{f}_n(\lambda)$
for a class of nonstationary nonlinear processes. Simple calculations
show that $\E\tilde{f}_n(\lambda)=\int_{0}^1f(t,\lambda) \,dt+o(1)$.
Hence $\tilde{f}_n(\lambda)$ is a consistent estimator of the averaged
energy at frequency $\lambda$ over time. An important observation from
Proposition \ref{prop:spectral-density} is that the asymptotic variance
of $\tilde{f}_n(\lambda)$ is determined by $\int_{0}^1f^2(t,\lambda
)
\,dt$, the averaged squared spectral density over time. The latter
quantity should be estimated if one wishes to construct confidence
intervals for $\int_{0}^1f(t,\lambda) \,dt$.


\section{Proofs}\label{sec:proof}
\mbox{}
\begin{pf*}{Proof of Theorem \ref{thm:clt_nondegen}}
Let $Z_k=\sum_{j=1}^nW_n(t_k,t_j)\{H_j(X_k)-\E[H_j(X_k)]\}$,
$k=1,2,\ldots, n$. Note that $N_n=\sum_{j=1}^nZ_j$. To prove the
theorem, we need to deal with the dependence structure of $\{Z_j\}$
first. According to \eqref{eq:fourier_repre},
\begin{eqnarray*}
Z_k&=&\int_{\R^2}g(t,s)\Biggl\{\sum
_{j=1}^nW_n(t_k,t_j)
\E\bigl[\beta_j(s)\bigr]\Biggr\} \gamma _k(t) \,dt\,ds\\[-2pt]
&:=&\int
_{\R^2}g(t,s)\Xi_{k}(s)\gamma_k(t) \,dt\,ds.
\end{eqnarray*}
Let $Z_{k,r}=\sum_{j=1}^nW_n(t_k,t_j)\{H_j(X_{k,r})-\E[H_j(X_{k,r})]\}
$. By Lemmas 1 and 2 in Zhou (\citeyear{Zho13}), the dependence measures
\begin{eqnarray*}
&&\|Z_{k}-Z_{k,r}\|_p\\[-2pt]
&&\qquad\le \int
_{\R^2}\bigl|g(t,s)\bigr|\bigl|\Xi_{k}(s)\bigr|\bigl\|\gamma
_k(t)-\gamma_{k,r}(t)\bigr\|_p \,dt\,ds
\\
&&\qquad\le \int
_{\R^2}\bigl|g(t,s)\bigr|\bigl|\Xi_{k}(s)\bigr|\bigl\{
\bigl\|L(X_k)\bigr\|_{2p}|t|^{\eta
}\bigl[\delta
_{X}(r,\eta2p)\bigr]^\eta+\delta_{L(X)}(r,p)\bigr
\} \,dt\,ds
\\
&&\qquad\le\int_{\R^2}C\bigl(1+|t|^{\eta}
\bigr)\bigl|g(t,s)\bigr|\bigl|\Xi_{k}(s)\bigr|\rho_1^r
\end{eqnarray*}
for some $\rho_1\in(0,1)$, where $p=2+\epsilon$. On the other hand,
note that
\[
\Xi_{k}(s)\le\sum_{j=1}^n\bigl|W_n(t_k,t_j)\bigr|\bigl|
\E\bigl[\beta_j(s)\bigr]\bigr|\le\sum_{j=1}^n\bigl|W_n(t_k,t_j)\bigr|
\bigl\|L(X_j)\bigr\|\le CW_{k,\cdot}.
\]
Hence $\|Z_{k}-Z_{k,r}\|_p\le CW_{k,\cdot}\rho_1^r$. As a second step,
we shall approximate $N_n$ by the sum of an $m$-dependent sequence. Define
$Z_{k,\{m\}}=\E[Z_k|\tilde{\FF}_{k,k-m}]$, where $\tilde{\FF
}_{k,k-m}=(\varepsilon_{k},\varepsilon_{k-1},\ldots,\varepsilon_{k-m})$.
For $j\in\Z$, define the projection operator
\[
\PP_{j}(\cdot)=\E[\cdot|\FF_j]-\E[\cdot|
\FF_{j-1}].
\]
Elementary manipulations show that $\PP_{k-r}Z_{k,\{m\}}=\E[\PP
_{k-r}Z_k|\tilde{\FF}_{k,k-m}]$. Hence by Jensen's inequality,
%
\begin{eqnarray}\quad
\label{eq:8} \bigl\|\PP_{k-r}[Z_k-Z_{k,\{m\}}]
\bigr\|_p&\le&\|\PP_{k-r}Z_k\|_p+\|\PP
_{k-r}Z_{k,\{m\}}\|_p\le2\|\PP_{k-r}Z_k
\|_p
\nonumber
\\[-8pt]
\\[-8pt]
\nonumber
&\le& 2\|Z_{k}-Z_{k,r}\|_p
\le CW_{k,\cdot}\rho_1^r.
\end{eqnarray}
Note that $Z_{k,\{m\}}-Z_k=\sum_{j=m}^\infty\{\E[Z_k|\tilde{\FF
}_{k,k-j}]-\E[Z_k|\tilde{\FF}_{k,k-j-1}]\}$
and the summands form a martingale difference sequence. By Burkholder's
inequality,
\begin{eqnarray*}
\|Z_{k,\{m\}}-Z_k\|^2_p&\le& C\sum
_{j=m}^\infty\bigl\|\E[Z_k|\tilde{\FF
}_{k,k-j}]-\E[Z_k|\tilde{\FF}_{k,k-j-1}]
\bigr\|_p^2
\\
&\le&C\sum_{j=m}^\infty
\|Z_k-Z_{k,j}\|_p^2\le
CW^2_{k,\cdot}\rho_1^{2m}.
\end{eqnarray*}
Therefore
%
\begin{equation}
\label{eq:9} \bigl\|\PP_{k-r}[Z_k-Z_{k,\{m\}}]
\bigr\|_p\le\|Z_{k,\{m\}}-Z_k\|_p\le
CW_{k,\cdot
}\rho_1^{m}.
\end{equation}
By \eqref{eq:8}, \eqref{eq:9} and Burkholder's inequality, for any
$r\ge0$,
\[
\Biggl\|\sum_{k=1}^n\PP_{k-r}[Z_k-Z_{k,\{m\}}]
\Biggr\|_p^2\le\sum_{k=1}^n
\Biggl\|\PP _{k-r}[Z_k-Z_{k,\{m\}}]\Biggr\|_p^2
\le CW^{(n)}\rho_1^{2\max(m,r)}.
\]
Observe that $\sum_{k=1}^n[Z_k-Z_{k,\{m\}}]=\sum_{r=0}^\infty\sum_{k=1}^n\PP_{k-r}[Z_k-Z_{k,\{m\}}]$. Therefore
%
\begin{eqnarray}
\label{eq:10}\Biggl \|\sum_{k=1}^n[Z_k-Z_{k,\{m\}}]
\Biggr\|_p&\le& \sum_{r=0}^\infty\Biggl\|\sum
_{k=1}^n\PP_{k-r}[Z_k-Z_{k,\{m\}}]
\Biggr\|_p
\nonumber
\\[-8pt]
\\[-8pt]
\nonumber
&\le& C\sqrt{W^{(n)}}m\rho_1^m.
\end{eqnarray}
Inequality \eqref{eq:10} shows that $N_n$ can be well approximated by
the sum of the $m$-dependent sequence $\{Z_{k,\{m\}}-\E[Z_{k,\{m\}}]\}
$. In particular, let $m=C\log n$. Then clearly approximation error in
\eqref{eq:10} can be made as $O(n^{-p})\sqrt{W^{(n)}}$ for any $p>0$.
In the final step we shall prove a central limit theorem for $N_{n,\{m\}
}:=\sum_{k=1}^n\{Z_{k,\{m\}}-\E[Z_{k,\{m\}}]\}$. Define the big blocks
and small blocks
\begin{eqnarray*}
R_j&=&\sum_{k=1}^{l_n}\bigl
\{Z_{(j-1)m+k,\{m\}}-\E[Z_{(j-1)m+k,\{m\}}]\bigr\} \quad\mbox {and }
\\
r_j&=&\sum_{k=l_n+1}^{s_n}\bigl
\{Z_{(j-1)m+k,\{m\}}-\E[Z_{(j-1)m+k,\{m\}
}]\bigr\},
\end{eqnarray*}
$j=1,2,\ldots,\lceil n/s_n\rceil$. Note that $R_j$'s are independent
and $r_j$'s are also independent. Then similar to the proof of \eqref
{eq:10}, we can obtain $|\var[N_n]-\var[N_{n,\{m\}}]|=o(W^{(n)})$, $\|
R_j\|^2_p=O(A_j)$,
%
\begin{equation}
\label{eq:11} \biggl\|\sum_jR_j
\biggr\|^2_p=O\biggl(\sum_jA_j
\biggr) \quad\mbox{and}\quad \biggl\|\sum_{j}r_j\biggr\|
^2_p=O\biggl(\sum_{j}a_j
\biggr).
\end{equation}
Therefore by condition (A3),
%
\begin{equation}
\label{eq:12} \biggl\|N_{n,\{m\}}-\sum_jR_j
\biggr\|_p=\biggl \|\sum_jr_j
\biggr\|_p=o\bigl(\sqrt{W^{(n)}}\bigr).
\end{equation}
By \eqref{eq:10}, $\|N_{n,\{m\}}-N_{n}\|_{p}=o(\sqrt{W^{(n)}})$. By
condition (A4), \eqref{eq:11} and \eqref{eq:12},
\[
\var\biggl[\sum_{j}R_j
\biggr]/W^{(n)}\ge c/2\qquad \mbox{for sufficiently large } n.
\]
Now by \eqref{eq:11},
\[
\frac{\sum_j\|R_j\|^p_p}{(\var[\sum_{j}R_j])^{p/2}}\le
C\frac{\sum_{j}A_j^{p/2}}{[W^{(n)}]^{p/2}}\le \biggl\{\frac{\max
A_j}{[W^{(n)}]} \biggr
\}^{p/2-1}\frac{\sum_{j}A_j}{W^{(n)}}.
\]
Hence by condition (A3), $\sum_j\|R_j\|^p_p/\var[\sum_{j}R_j])^{p/2}\rightarrow0$.
Therefore by the Lyapunov CLT, $\sum_{j}R_j/\sqrt{\var[\sum_{j}R_j]}\Rightarrow N(0,1)$. Now by \eqref{eq:10}--\eqref{eq:12}, the
theorem follows.
\end{pf*}

\begin{pf*}{Proof of Theorem \ref{thm:order_degen}}
For any fixed $(x,y)\in\R^2$, define
\[
\varrho_n(x,y)=\sum_{k,j}W_n(t_k,t_j)
\gamma_k(x)\gamma_j(y).
\]
We shall first determine the order of magnitude of $\varrho_n(x,y)$.
For complex-valued random variables $X$ and $Y$, define $V(X)=\E|X-\E
X|^2$ and\break  $\cov[X,Y]=\E[(X-\E X)(\bar{Y}-\E\bar{Y})]$. Note that, by
(A6) and the symmetry of $W(\cdot,\cdot)$,
\begin{eqnarray*}
\label{eq:13} V\bigl(\varrho_n(x,y)\bigr)&=&\sum
_{k,j}\sum_{l,m}W_n(t_l,t_m)W_n(t_k,t_j)
\cov \bigl[\gamma_l(x)\gamma_m(y),\gamma_k(x)
\gamma_j(y)\bigr]
\\
&\le&\sum_{k,j}\sum
_{l,m}W^2_n(t_k,t_j)\bigl|
\cov\bigl[\gamma_l(x)\gamma _m(y),\gamma_k(x)
\gamma_j(y)\bigr]\bigr|
\\
&&{}+2^{p/2}\sum
_{k,j}\sum_{l,m}f_n(t_k,t_j)\bigl|W_n(t_k,t_j)\bigr|
\rho ^p_{k,j}(l,m)\\
&&\hspace*{63pt}{}\times\bigl|\cov\bigl[\gamma_l(x)
\gamma_m(y),\gamma_k(x)\gamma _j(y)\bigr]\bigr|
\\
&:=&\sum_{k,j}W^2_n(t_k,t_j)
\varrho_{k,j}(x,y;0)\\
&&{}+2^{p/2}\sum_{k,j}f_n(t_k,t_j)\bigl|W_n(t_k,t_j)\bigr|
\varrho_{k,j}(x,y;p),
\end{eqnarray*}
where $\rho_{k,j}(l,m)=\min \{\max(|l-k|,|m-j|),\max
(|l-j|,|m-k|)
\}$ and
\[
\varrho_{k,j}(x,y;p)=\sum_{l,m}
\rho^p_{k,j}(l,m)\bigl|\cov\bigl[\gamma _l(x)\gamma
_m(y),\gamma_k(x)\gamma_j(y)\bigr]\bigr|.
\]
We will omit the subscripts $k,j$ in $\rho_{k,j}(l,m)$ in the sequel
for simplicity. Let
\begin{eqnarray*}
&&\rho^*(l,m)=\max \bigl\{\min\bigl(|k-l|,|k-m|,|j-l|,|j-m|\bigr),\\
&&\hspace*{75pt}\min
\bigl(|l-k|,|l-j|,|l-m|\bigr),
\min\bigl(|m-k|,|m-j|,|m-l|\bigr),\\
&&\hspace*{79pt}\min\bigl(|k-l|,|k-m|,|k-j|\bigr),\min
\bigl(|j-l|,|j-m|,|j-k|\bigr) \bigr\}.
\end{eqnarray*}
We shall first show that
%
\begin{equation}
\label{eq:14} \rho(l,m)\le2\rho^*(l,m) \qquad\mbox{for all } l,m.
\end{equation}
By the symmetry of $\rho(l,m)$ and $\rho^*(l,m)$, we only need to
consider the case $l\le m$ and $k \le j$. Now if $l\le k$, then
\begin{eqnarray*}
2\min\bigl(|l-k|,|l-j|,|l-m|\bigr)&\ge&|k-l| \qquad\mbox{if }
 m\ge(l+k)/2 \quad\mbox{ and}
\\
2\min\bigl(|k-l|,|k-m|,|j-l|,|j-m|\bigr)&\ge&|l-k| \qquad\mbox{if } m\le(l+k)/2.
\end{eqnarray*}
If $l\ge k$, then
\begin{eqnarray*}
2\min\bigl(|k-l|,|k-m|,|k-j|\bigr)&\ge&|k-l| \qquad \mbox{if } j\ge(l+k)/2 \quad\mbox{ and }
\\
 2\min\bigl(|k-l|,|k-m|,|j-l|,|j-m|\bigr)&\ge&|l-k| \qquad\mbox{if } j\le(l+k)/2.
\end{eqnarray*}
In summary, $2\rho^*(l,m)\ge|l-k|$. Similarly, $2\rho^*(l,m)\ge
|j-m|$. Hence \eqref{eq:14} follows. Now by Lemmas 1, 2 and 7 in Zhou
(\citeyear{Zho13}), elementary calculations using condition (A5) show that
%
\begin{eqnarray}
\label{eq:15}\bigl |\cov\bigl[\gamma_l(x)\gamma_m(y),
\gamma_k(x)\gamma_j(y)\bigr]\bigr|&\le& C\bigl(|x|^{\eta
}+1
\bigr) \bigl(|y|^{\eta}+1\bigr)r_1^{\rho^*(l,m)}
\nonumber
\\[-8pt]
\\[-8pt]
\nonumber
&\le& C
\bigl(\bigl|(x,y)\bigr|^{2\eta}+1\bigr)r_1^{\rho^*(l,m)}
\end{eqnarray}
for some $r_1\in[0,1)$. Observe that for $r=0,1,\ldots,n$, the number
of pairs $(l,m)$ such that $\rho(l,m)=r$ is at most $2r^2$. Now by
\eqref{eq:14} and \eqref{eq:15}, we obtain that, for $r_2=\sqrt{r_1}$,
%
\begin{eqnarray}
\label{eq:16} \varrho_{k,j}(x,y;p)&\le& C\sum
_{l,m}\rho^p(l,m) \bigl(\bigl|(x,y)\bigr|^{2\eta
}+1
\bigr)r_2^{\rho(l,m)}
\nonumber
\\[-8pt]
\\[-8pt]
\nonumber
&\le& C\bigl(\bigl|(x,y)\bigr|^{2\eta}+1
\bigr).
\end{eqnarray}
Note that the constant $C$ does not depend on $(k,j)$. Now by \eqref
{eq:16}, we obtain
\[
\bigl\|\varrho_n(x,y)-\E\varrho_n(x,y)\bigr\|\le C
\bigl(\bigl|(x,y)\bigr|^{\eta}+1\bigr)\sqrt {W_{(n)}+\Delta_n}.
\]
Therefore by condition (A1),
\begin{eqnarray*}
\|D_n-\E D_n\|&\le&\int_{\R^2}\bigl|g(x,y)\bigr|
\bigl\|\varrho_n(x,y)-\E\varrho _n(x,y)\bigr\| \,dx \,dy
\\
&\le& C
\sqrt{W_{(n)}+\Delta_n}\int_{\R^2}
\bigl(\bigl|(x,y)\bigr|^{\eta
}+1\bigr)\bigl|g(x,y)\bigr| \,dx \,dy\\
&\le& C\sqrt{W_{(n)}+
\Delta_n}.
\end{eqnarray*}
The theorem follows.
\end{pf*}

\begin{pf*}{Proof of Theorem \ref{thm:degen}}
Let $\varpi_n(x,y)=\sum_{k,j=1}^nW_n(t_k,t_j)\Gamma_n^*(t_k,x)\Gamma
_n^*(t_j,y)$. Recall the definition of $\Gamma_n^*(t,s)$ in \eqref
{eq:40}. We shall show that the two processes $\varrho_{n}(x,y)$ and
$\varpi_n(x,y)$ are close in the sense that
%
\begin{eqnarray}
\label{eq:23}&&\biggl |\E h\biggl[\int_{\R^2}g(x,y)\bigl[
\varrho_n(x,y)-\E\varrho_n(x,y)\bigr] \,dx\,dy\biggr]
\nonumber
\\[-8pt]
\\[-8pt]
\nonumber
&&\qquad{}-
\E h\biggl[\int_{\R^2}g(x,y)\bigl[\varpi_n(x,y)-\E
\varpi_n(x,y)\bigr] \,dx\,dy\biggr] \biggr|\rightarrow0
\end{eqnarray}
for any bounded and continuous $h$. For any $s>0$, define the region
$A(s)=\{(x,y)\in\R^2, |x|\le s, |y|\le s\}$, and let $\bar{A}(s)=\R
^2/A(s)$. Note that, by Lemma 4 in Zhou (\citeyear{Zho13}), we have
\begin{eqnarray*}
&&\biggl\|\int_{\bar{A}(s)}g(x,y)\bigl[\varrho_n(x,y)-\E
\varrho_n(x,y)\bigr] \,dx\,dy\biggr\| \\
&&\qquad\le \int_{\bar{A}(s)}\bigl|g(x,y)\bigr|
\bigl\|\varrho_n(x,y)-\E\varrho_n(x,y)\bigr\| \,dx\,dy
\\
&&\qquad\le C\int
_{\bar{A}(s)}\bigl|g(x,y)\bigr|\bigl(1+\bigl|(x,y)\bigr|^{\eta}\bigr) \,dx\,dy.
\end{eqnarray*}
Observe that $\int_{\bar{A}(s)}|g(x,y)|(1+|(x,y)|^{\eta}) \,dx\,dy$ is
independent of $n$ and converges to $0$ as $s\rightarrow\infty$.
Similar inequality holds for $\int_{\bar{A}(s)}g(x,y)[\varpi
_n(x,y)-\E
\varpi_n(x,y)] \,dx\,dy$. Hence, to prove \eqref{eq:23}, one only need to
show that, for each fixed $s$, \eqref{eq:23} holds with $\R^2$ therein
replaced by $A(s)$. To this end, we shall first show that, for any
$(x_1,y_1),\ldots,(x_m,y_m)\in A(s)$ and any bounded and continuous $h\dvtx \R^{m}\rightarrow\R$,
%
\begin{eqnarray}
\label{eq:24}\quad  &&\E h\bigl\{\bigl(\varrho_n(x_1,y_1)-
\E\varrho_n(x_1,y_1),\ldots,\varrho
_n(x_m,y_m)-\E\varrho_n(x_m,y_m)
\bigr)\bigr\}\nonumber\\
&&\quad{}-\E h\bigl\{\bigl(\varpi_n(x_1,y_1)-
\E\varpi_n(x_1,y_1),\ldots,\varpi
_n(x_m,y_m)-\E\varpi_n(x_m,y_m)
\bigr)\bigr\}\\
&&\qquad=o(1).\nonumber
\end{eqnarray}
We shall only prove the case $m=1$ since similar arguments apply to
general $m$. Consider the case $x,y\neq0$. By Corollary 2 of \citet{WuZho11}, we have, on a possibly richer probability space, a
sequence a i.i.d. 4-dimensional standard normal random vectors $\ZZ
_1,\ZZ_2,\ldots, \ZZ_n$, such that
%
\begin{eqnarray}
\label{eq:30} &&\max_{1\le k\le n} \Biggl|\sum
_{j=1}^k \bigl\{\bigl(\gamma_{j,1}(x),
\gamma _{j,2}(x),\gamma_{j,1}(y),\gamma_{j,2}(y)
\bigr)^\top
-\Sigma^{1/2}_{j^*}
\bigl(t_j,(x,y)^\top\bigr)\ZZ_j \bigr\}\Biggr|
\nonumber
\\[-8pt]
\\[-8pt]
\nonumber
&&\qquad=o_\p \bigl(n^{1/4}\log ^2 n\bigr).
\end{eqnarray}
Define the complex-valued random variables $Z^*_j=[\Sigma
^{1/2}_{j^*}(t_j,(x,y)^\top)\ZZ_j]_{1}+i[\Sigma
^{1/2}_{j^*}(t_j,(x,y)^\top)\ZZ_j]_{2}$ and $Z^{**}_j=[\Sigma
^{1/2}_{j^*}(t_j,(x,y)^\top)\ZZ_j]_{3}+i[\Sigma
^{1/2}_{j^*}(t_j,(x,y)^\top)\times \ZZ_j]_{4}$, where $[\Sigma
^{1/2}_{j^*}(t_j,(x,y)^\top)\ZZ_j]_{r}$ denotes the $r$th element of
$\Sigma^{1/2}_{j^*}(t_j,(x,y)^\top)\times \ZZ_j$.
Define the quadratic form
%
\begin{equation}
\label{eq:36} \varrho_n^{\dmd}(x,y)=\sum
_{k,j=1}^nW_n(t_k,t_j)Z^*_kZ^{**}_j.
\end{equation}
Note that
\begin{eqnarray*}
\bigl|\varrho_n(x,y)-\varrho_n^{\dmd}(x,y)\bigr|&\le&\biggl|\sum
_{k,j}W_n(t_k,t_j)
\bigl[\gamma _{k}(x)-Z^*_{k}\bigr]\gamma_{j}(y)\biggr|\\
&&{}+\biggl|
\sum_{k,j}W_n(t_k,t_j)
\bigl[\gamma _{j}(y)-Z^{**}_{j}
\bigr]Z^{*}_{k}\biggr|.
\end{eqnarray*}
Using the summation by parts technique and similar to the proof of
inequality~(12) in Zhou (\citeyear{Zho13}), we have $|\varrho_n(x,y)-\varrho
_n^{\dmd}(x,y)|=o_\p(1)$.
Now by the above inequality and Lemma 3 in Zhou (\citeyear{Zho13}), we obtain
%
\begin{equation}
\label{eq:31} \bigl|\varrho_n(x,y)-\E\varrho_n(x,y)-\bigl[
\varrho_n^{\dmd}(x,y)-\E\varrho _n^{\dmd}(x,y)
\bigr]\bigr|=o_\p(1).
\end{equation}
By condition (A10), we have that $\Sigma_k(t,(x,y)^{\top})$ is
continuous on $[b_{k},b_{k+1}]$, $k=0,1,\ldots,r$. Therefore elementary
calculations show that, on a possibly richer probability space,
%
\begin{equation}
\label{eq:32} \bigl|\varrho_n^{\dmd}(x,y)-\E\varrho_n^{\dmd}(x,y)-
\bigl[\varpi_n(x,y)-\E \varpi _n(x,y)\bigr]\bigr|=o_\p(1).
\end{equation}
Hence \eqref{eq:24} with $m=1$ follows from \eqref{eq:31} and \eqref
{eq:32}. Now consider the case $x=0$. If $L(\cdot)\equiv C$, then
$\varrho_n(0,y)\equiv0$ and $\varpi_n(0,y)\equiv0$. Hence \eqref
{eq:24} trivially holds. If $L(\cdot)$ is not a constant, then \eqref
{eq:24} follows from similar and simpler arguments as above by
considering the covariance matrix $\Sigma^*_k(y)$. In summary, \eqref
{eq:24} follows.

As a second step toward \eqref{eq:23}, we prove that $\{\varrho
_n(x,y)-\E\varrho_n(x,y)\}$ is tight on ${\cal C}(A(s))$, where
${\cal
C}(A(s))$ is the collection of all complex-valued continuous functions
on $A(s)$ equipped with the uniform topology. Note that
%
\begin{eqnarray}
\label{eq:33} &&\varrho_n(x_1,y_1)-\E
\varrho_n(x_1,y_1)-\bigl[
\varrho_n(x_2,y_2)-\E \varrho_n(x_2,y_2)
\bigr]
\nonumber
\\
&&\qquad=i\int_{y_2}^{y_1}\rho^{(1)}_n(x_2,y)
\,dy
+i\int_{x_2}^{x_1}\rho^{(2)}_n(x,y_2)
\,dx\\
&&\qquad\quad{}-\int_{x_2}^{x_1}\int_{y_2}^{y_1}
\rho^{(3)}_n(x,y) \,dx\,dy,\nonumber
\end{eqnarray}
where $\rho^{(1)}_n(x,y)=\sum_{k,j}W_n(t_k,t_j)[\gamma_k(x)\gamma
^{\dmd
}_j(y)-\E\gamma_k(x)\gamma^{\dmd}_j(y)]$, $\rho^{(2)}_n(x,y)=\sum_{k,j}W_n(t_k,t_j)\times$ $[\gamma^{\dmd}_k(x)\gamma_j(y)-\E\gamma
^{\dmd}_k(x)\gamma_j(y)]$ and
\[
\rho^{(3)}_n(x,y)=\sum_{k,j}W_n(t_k,t_j)
\bigl[\gamma^{\dmd}_k(x)\gamma ^{\dmd
}_j(y)-
\E\gamma^{\dmd}_k(x)\gamma^{\dmd}_j(y)
\bigr]
\]
with $\gamma^{\dmd}_j(x)=L^*(X_j)e^{ixX_j}-\E[L^*(X_j)e^{ixX_j}]$. By
the proof of Lemma 4 in Zhou (\citeyear{Zho13}), we have
%
\begin{equation}
\label{eq:34} \sup_{k=1,2,3, (x,y)\in A(s)}\bigl\|\rho^{(k)}_n(x,y)
\bigr\|_{2+\epsilon}=O(1).
\end{equation}
By \eqref{eq:33} and \eqref{eq:34}, we have, for any fixed
$(x_0,y_0)\in A(s)$ and $\delta>0$,
%
\begin{eqnarray}
\label{eq:35} &&\Bigl \|\sup_{|(x,y)-(x_0,y_0)|\le\delta}\bigl|\varrho_n(x,y)-\E
\varrho _n(x,y)
-\bigl[\varrho_n(x_0,y_0)-
\E\varrho_n(x_0,y_0)\bigr]\bigr|
\Bigr\|_{2+\epsilon
}
\nonumber
\\[-8pt]
\\[-8pt]
\nonumber
&&\qquad=O(\delta).
\end{eqnarray}
Define $\omega(\delta)=\sup_{|(x_1,y_1),(x_2,y_2)\in A(s),
|(x_1,y_1)-(x_2,y_2)|\le\delta|}|\varrho_n(x_1,y_1)-\E\varrho
_n(x_1,y_1)-[\varrho_n(x_2,y_2)-\E\varrho_n(x_2,y_2)]|$. By \eqref
{eq:35} and a standard chaining technique, we have for each fixed
$\varepsilon>0$
\[
\lim_{\delta\rightarrow0}\limsup_{n\rightarrow\infty} \p\bigl(\omega (
\delta )>\varepsilon\bigr)=0.
\]
Hence $\{\varrho_n(x,y)-\E\varrho_n(x,y)\}$ is tight on $A(s)$. By
standard smooth Gaussian process techniques, it is easy to see that
$\varpi_n(x,y)-\E\varpi_n(x,y)$ is tight on $A(s)$. Since both
processes are relatively compact on $A(s)$ and the differences of their
finite dimensional distributions converge in the sense of \eqref
{eq:24}, we have for any bounded and continuous function $h^*$,
\[
\E h^*\bigl(\varrho_n(x,y)-\E\varrho_n(x,y)\bigr)- \E
h^*\bigl(\varpi_n(x,y)-\E \varpi _n(x,y)\bigr)
\rightarrow0.
\]
Since $g(x,y)\in L^1(\R^2)$, we have $K(f):=\int_{A(s)}g(x,y)f(x,y)
\,dx\,dy$ is continuous on ${\cal C}(A(s))$. Hence \eqref{eq:23} follows.

Finally, note that $\Gamma_n^*(t_k,s)$'s are independent complex-valued
Gaussian processes, $k=0,1,2,\ldots,n$. By the classic Gaussian process
theory [see, e.g., \citet{Kuo75}, Chapter~1.2], we have $\Gamma
_n^*(t_k,s)$ can be represented as (in the sense of equality in distribution)
\[
\bigl(\Re\bigl(\Gamma_{n}^*(t_k,s)\bigr),\Im\bigl(
\Gamma_{n}^*(t_k,s)\bigr)\bigr)^\top=\sum
_{j=1}^\infty A_{n,j}(t_k,s)B_{k,j},\qquad
k=0,1,2,\ldots,n,
\]
where $\Re(\cdot)$ and $\Im(\cdot)$ denotes real and imaginary
parts of
a complex number, respectively, $A_{n,j}(t,s)$'s are $2\times2$ matrix
functions and $B_{k,j}$'s are independent 2-dimensional standard normal
random vectors. Hence it is straightforward to see that $\int_{\R
^2}g(x,y)[\varpi_n(x,y)-\E\varpi_n(x,y)] \,dx\,dy$ is a quadratic form of
i.i.d. standard normal random variables $G_1,G_2,\ldots.$ Moreover, by
the arguments of Lemma 4 in Zhou (\citeyear{Zho13}), we have $\|\int_{\R
^2}g(x,y)[\varpi_n(x,y)-\E\varpi_n(x,y)] \,dx\,dy\|=O(1)$. Hence
\[
\int_{\R^2}g(x,y)\bigl[\varpi_n(x,y)-\E
\varpi_n(x,y)\bigr] \,dx\,dy=\sum_{j=1}^\infty
\alpha_{n,j}\bigl(Z_j^2-1\bigr)
\]
with $\sum_{j=1}^\infty\alpha^2_{n,j}<\infty$.
\end{pf*}
%
\begin{remark}
As we can see from the proof of \eqref{eq:30}, the
positive-definiteness requirement on $\Sigma_k(t,s)$ and $\Sigma
^*_k(t,s)$ in (A9) is to facilitate a Gaussian approximation result in
\citet{WuZho11}. We point out that the positive-definiteness
requirement can be weakened to the assumption that certain block sums
of the latter long-run covariance matrices are positive definite. See
Remark 2 of \citet{WuZho11}. For presentational simplicity, we shall
stick to the everywhere positive definiteness assumption in this paper.
\end{remark}

\begin{pf*}{Proof of Corollary \ref{coro:clt_degen}}
We shall prove this corollary by showing that $\varrho_n(x,y)-\E
\varrho
_n(x,y)$ converges to a Gaussian measure\vadjust{\goodbreak} on ${\cal C}(A(s))$. By the
tightness of $\{\varrho_n(x,y)-\E\varrho_n(x,y)\}$ and the arguments
in the proof of Theorem \ref{thm:degen}, it suffices to show that any
finite dimensional distribution of the latter sequence of measures
converges to a (multivariate) normal distribution. To this end, we will
only show that $\rho^*_n(x,y)-\E\rho^*_n(x,y)$ converges to a Gaussian
distribution for any $(x,y)\in A(s)$ since all other cases follow by
similar arguments and the Cramer--Wold device. Here $\rho
^*_n(x,y)=\sum_{k,j=1}^nW_n(t_k,t_j)\gamma_{k,1}(x)\gamma_{j,1}(y)$, $\gamma
_{k,1}(x)=\Re(\gamma_{k}(x))$ and $\gamma_{k,2}(x)=\Im(\gamma_{k}(x))$.
Consider the case $x,y\neq0$. Then by the proof of \eqref{eq:31}, we have
\[
\bigl|\rho^*_n(x,y)-\E\rho^*_n(x,y)-\bigl[
\rho_{n,1}^{\dmd}(x,y)-\E\rho _{n,1}^{\dmd}(x,y)
\bigr]\bigr|=o_\p(1),
\]
where $\rho_{n,1}^{\dmd}(x,y)=\sum_{k,j=1}^nW_n(t_k,t_j)\Re
(Z^*_k)\Re
(Z^{**}_j)$. Recall the definitions of $Z^*_k$ and $Z^{**}_j$ in \eqref
{eq:36}. Note that $\rho_{n,1}^{\dmd}(x,y)$ is a quadratic form of
i.i.d. normal random variables. More specifically, $\rho_{n,1}^{\dmd
}(x,y)$ can be written as
%
\begin{equation}
\label{eq:44} 2\rho_{n,1}^{\dmd}(x,y)=\ZZ^\top
D^\top(W_n\otimes A) D \ZZ,
\end{equation}
where $\ZZ=(\ZZ^\top_1,\ldots,\ZZ^\top_n)^\top$ is a length $4n$
vector of
i.i.d. standard normal random variables, $D=\operatorname{Diag}(\Sigma
^{1/2}_{1^*}(t_1,(x,y)^\top),\ldots, \Sigma
^{1/2}_{n^*}(t_n,(x,y)^\top
))$ is a $4n\times4n$ block diagonal matrix, $A$ is the $4\times4$
matrix with $(1,3)$th and $(3,1)$th elements equaling 1 and all other
entries equaling 0 and $\otimes$ denotes the Kronecker product. Let
$M_n=D^\top W_n\otimes A D$. By condition (A9), $D$ is positive
definite with eigenvalues bounded both above and below. Then it is easy
to see that there exist constants $0<c\le C<\infty$, such that $c\le
\sum_{k,j=1}^{4n}M^2_n(k,j)\le C$. By the Lyapunov CLT, to prove the
asymptotic normality of $\rho_{n,1}^{\dmd}(x,y)-\E\rho_{n,1}^{\dmd
}(x,y)$, it suffices to show that the $|\theta^*_{n,1}|\rightarrow0$,
where $\theta_{n,1}^*$ is the eigenvalue of $W_n\otimes A$ with the
maximum absolute value. By the basic property of Kronecker product, we
have that the eigenvalues of $W_n\otimes A$ are the products of the
eigenvalues of $W_n$ and $A$. Hence it is clear that $|\theta
^*_{n,1}|=|\theta_{n,1}|\rightarrow0$. The case when $x=0$ or $y=0$
follows similarly.
\end{pf*}


\begin{pf*}{Sketch of Proof of Theorem \ref{thm:ftransform_normal_appr}}
Theorem \ref{thm:ftransform_normal_appr} follows from Lemma 5 in the online supplement of the paper with $|b-a|=\lf c^{1/2}n^{1/2}\log^{-1}
n \rf$ for some finite constant $c$ together with a careful check of
the proof of Theorem 1 in \citet{WuZho11} with $l=\lf c\log n\rf$
and $m=\lf l^{1/2}n^{1/2}\log^{-3/2} n \rf$ therein.
\end{pf*}


\section*{Acknowledgments} The author is grateful to the two
anonymous referees for their careful reading of the manuscript and many
helpful comments.

\begin{supplement}[id=suppA]
\stitle{Supplement for ``Inference of weighted $V$-statistics for
nonstationary time series and its applications''}
\slink[doi]{10.1214/13-AOS1184SUPP} 
\sdatatype{.pdf}
\sfilename{aos1184\_supp.pdf}
\sdescription{This supplementary material contains auxiliary lemmas
and proofs of Propositions \ref{prop:ac}, \ref{prop:local-linear},
\ref{prop:spectral-density} and Corollaries \ref{coro:nonclt-degen1},
\ref{coro:nonclt-degen2} of the paper.}
\end{supplement}


%

\printaddresses

\end{document}